\documentclass[11pt]{article}

\usepackage{amsmath,amssymb,amsthm}
\usepackage[hidelinks]{hyperref}
\usepackage{geometry}
\date{}

\newcommand{\N}{\mathbb{N}}

\newcommand{\R}{\mathbb{R}}
\newcommand{\Z}{\mathbb{Z}}

\newcommand{\be}{\begin{equation}}
\newcommand{\ee}{\end{equation}}

\newcommand{\K}{{\cal K}}

\newcommand{\T}{\mathbb{T}}

\def\g{{\gamma }}

\def\){\right)}
\def\({\left(}

\def\a{{\alpha}}
\def\b{{\beta}}

\def\vp{{\varphi}}

\def\g{{\gamma }}

\def\A{{\mathcal{A}}}
\def\I{{\rm{I}}}

\numberwithin{equation}{section}

\newtheorem{definition}{Definition}[section]
\newtheorem{remark}[definition]{Remark}
\newtheorem{theorem}[definition]{Theorem}
\newtheorem{lemma}[definition]{Lemma}

\newtheorem{corollary}[definition]{Corollary}

\title{Sharp error estimates for Lagrange interpolation
and Kantorovich-type operators in weighted $L_p$-spaces}

\author{
\textbf{L. Angeloni$^{1}$, D. Costarelli$^{1}$,
Yu. Kolomoitsev$^{2}$\thanks{Corresponding author.}}
\\[2ex]
$^{1}$Department of Mathematics and Computer Science\\
University of Perugia\\
Via Vanvitelli 1, 06123 Perugia, Italy\\
{\small\ttfamily
laura.angeloni@unipg.it, danilo.costarelli@unipg.it}
\\[2ex]
$^{2}$Institute for Numerical and Applied Mathematics\\
University of G\"ottingen\\
Lotzestr. 16--18, 37083 G\"ottingen, Germany\\
{\small\ttfamily kolomoitsev@math.uni-goettingen.de}
}

\date{}

\begin{document}

\maketitle
\begin{abstract}
In this paper, we establish sharp direct and inverse approximation estimates for algebraic Lagrange interpolation polynomials and Kantorovich-type operators in weighted $L_p$ spaces on $[-1,1]$.
 A well-known limitation in the theory of approximation by interpolation polynomials states that for small rates of convergence $\alpha \le 1/p$, classical moduli of smoothness fail to provide matching direct and inverse estimates due to the insufficiency of integral metrics to determine pointwise values at the interpolation nodes. To overcome this difficulty, we introduce new special semi-discrete measures of smoothness based on Steklov averaging operators,
which combine global $L_p$-smoothness
with discrete deviations from Steklov
averages at the interpolation nodes.
 By using these measures, we establish fully matching direct and inverse approximation theorems for Lagrange interpolation polynomials, thereby providing complete two-sided characterizations of the rate of convergence for an arbitrary rate $\alpha > 0$. As a natural application of the developed approach, we also obtain corresponding direct and inverse estimates, as well as strong converse inequalities, for Kantorovich-type operators constructed via Lagrange polynomials.

\end{abstract}

\section{Introduction} \label{sec1}

In the constructive theory of functions, establishing a complete characterization of the rate of convergence for a given approximation process via matching direct (Jackson-type) and inverse (Bernstein-type) estimates is a fundamental objective. For the best algebraic polynomial approximation, such estimates have long been established in both unweighted and weighted $L_p[-1,1]$ spaces. Typically, these results are formulated in terms of the Ditzian--Totik modulus of smoothness
$$
\omega_\varphi^r (f,t)_{p}:=\sup_{0 < h \le t} \| \Delta^r_{h\varphi} f \|_{L_p[-1,1]},\quad \varphi(x) := \sqrt{1 - x^2}.
$$
Unlike the classical moduli of smoothness, this modulus specifically captures the boundary behavior of functions near the endpoints of the interval, thereby providing sharp error estimates for algebraic polynomial approximation.

In this setting, we recall the following well-known inequalities:
\begin{equation}\label{intrjack}
    \|f-P_{n,f}\|_{L_p[-1,1]} \le C \omega_\varphi^r (f,1/n)_{p}
\end{equation}
and
\begin{equation}\label{intrinv}
    \omega_\varphi^r (f,1/n)_{p} \le \frac{C}{n^r}\sum_{k=0}^n (k+1)^{r-1}\|f-P_{k,f}\|_{L_p[-1,1]},
\end{equation}
where $P_{n,f} \in \mathcal{P}_n$ is the polynomial of best approximation to $f$ in $L_p[-1,1]$ and $C$ is a positive constant independent of $f$ and $n$ (see, e.g.,~\cite[Theorems~7.2.1 and~7.2.4]{DITO1987}; see also~\cite[Sec.~2.5.2]{MAMI2008} for the corresponding results in weighted $L_p$ spaces). As usual, estimates \eqref{intrjack} and \eqref{intrinv} imply the following standard equivalence for all $\alpha \in (0,r)$:
\begin{equation}\label{intrequiv}
    \|f-P_{n,f}\|_{L_p[-1,1]} = \mathcal{O}(n^{-\alpha}) \quad \iff \quad \omega_\varphi^r (f,t)_{p} = \mathcal{O}(t^\alpha).
\end{equation}

Similar direct and inverse estimates (and even strong converse inequalities) leading to characterizations of type~\eqref{intrequiv} for an arbitrary $\alpha > 0$ have also been obtained for various approximation methods such as linear summation methods of Fourier--Jacobi series (see, e.g.,~\cite{BA1974, DI1998, MATH2008}), certain Kantorovich-type operators~\cite{CHDI1994, GOZH1995}, and  Durrmeyer-type polynomial operators~\cite{CHDIIV1993}. Of course, for a given approximation method, one needs to employ a specific modulus of smoothness or $K$-functional that appropriately captures the boundary behavior of the function near the endpoints of the interval $[-1,1]$, in the spirit of the Ditzian--Totik moduli (see the survey~\cite{DI2007}).

In contrast, for interpolation and sampling operators, establishing fully matching direct and inverse estimates that guarantee an equivalence of type~\eqref{intrequiv} for an arbitrary order $\alpha > 0$ has not been completely achieved. Most existing results in this direction are essentially restricted to the case $\alpha > 1/p$; see, e.g.,~\cite{TH2018, IV1986, MARU1999}.
For example, let $\mathcal{L}_n(f,u)$ be the Lagrange
polynomial interpolating $f\in C[-1,1]$ at the zeros
of the $n$th orthogonal polynomial associated with $u$.
Under the assumptions on $u$ of
\cite[Theorem~3.1]{MARU1999} with $w=1$,
the direct estimate~\eqref{eq:MR1999}
and the inverse inequality~\eqref{intrinv}
yield, for $1<p<\infty$ and $\alpha\in(1/p,r)$,
\begin{equation}\label{intrequivL}
\|f-\mathcal{L}_n(f,u)\|_{L_p[-1,1]}
=\mathcal{O}(n^{-\alpha})
\quad\iff\quad
\omega_\varphi^r(f,t)_p=\mathcal{O}(t^\alpha).
\end{equation}

However, already at the critical threshold
$\alpha=1/p$, standard measures of smoothness
fail to yield the desired equivalence.
The known error estimates in this case have been
restricted to direct estimates for functions of
bounded variation; see, e.g.,~\cite{PR93,MARU2010}.
Related results for bivariate polynomial interpolation
and functions of bounded variation were obtained
in~\cite{KLP2018}.
 At the same time, it is well known (see~\cite{AS1973, NE1984}) that under certain natural restrictions on the weight $u$, the Lagrange polynomials $\mathcal{L}_n(f,u)$ converge in the $L_p$ norm for any continuous function $f$ on $[-1,1]$. Thus, for $\alpha \in (0,1/p]$, there is a significant gap in the theory concerning the complete characterization of the approximation properties of $\mathcal{L}_n(f,u)$.

One of the reasons for this gap lies in
the fundamental difference between
integral smoothness and pointwise information.
While the standard Ditzian--Totik moduli of
smoothness and other related functionals
measure the smoothness properties of a function
in an integral sense
(see, e.g.,~\cite{POKA1998,IV1982,BU1983,DITO1987}),
the Lagrange interpolation process depends
on the exact values of $f$ at discrete nodes.
This distinction becomes particularly important
in the low-smoothness range $\alpha\le 1/p$,
where integral $L_p$-smoothness alone does
not, in general, control these point values.

The main objective of the present paper is to extend  the characterization of convergence
rates for Lagrange interpolation to the range
$0<\alpha\le 1/p$, thereby obtaining a
characterization valid for all $\alpha>0$. To this end, we introduce a new class of
measures of smoothness based on a special Steklov averaging operator and weighted
discrete norms over the interpolation nodes. These measures are adapted
both to the weighted $L_p$ metric and to the discrete geometry of the
interpolation nodes. Our approach follows the ideas of~\cite{KL23} and~\cite{K25} (see also~\cite{KO2026}), where analogous problems were solved for the
approximation of functions on the torus $\mathbb{T}$ and on the real
line $\mathbb{R}$, respectively.
More precisely, using the Steklov averaging operator defined by
\begin{equation*}
\mathcal{A}_h f(x)
:=
\frac{1}{2h}
\int_{\arccos x-h}^{\arccos x+h}
f(\cos v)\,dv,
\qquad x\in[-1,1],
\end{equation*}
we obtain
matching direct and inverse approximation theorems for approximation by $\mathcal{L}_n(f,u)$ in $L_{p,w}[-1,1]$
(see Theorems~\ref{thsp} and~\ref{prinvp}). For simplicity, we present
here their unweighted form and omit the technical assumptions on the
interpolation weight $u$. In particular, for every $\alpha\in(0,2)$
and $p>2$, we obtain equivalence relations of the form
\begin{equation}\label{intrequivL2}
\|f-\mathcal{L}_n(f,u)\|_{L_p[-1,1]}
=
\mathcal{O}(n^{-\alpha})
\quad\Longleftrightarrow\quad
\omega_\varphi^2(f,1/n)_p
+
\|f-\mathcal{A}_{1/n}f\|_{\ell_p(X_n(u))}
=
\mathcal{O}(n^{-\alpha}),
\end{equation}
where $\|\cdot\|_{\ell_p(X_n(u))}$ denotes the discrete norm associated
with the interpolation nodes $X_n(u)$; see~\eqref{discrnormlp}.

Furthermore, as a natural application of the developed technique, we establish sharp direct and inverse estimates, as well as strong converse inequalities, for the corresponding Kantorovich-type operators (see Subsection~\ref{secKant}).

The paper is organized as follows.
Section~2 contains preliminary definitions and results.
We introduce the relevant weight classes and conditions
in Subsection~\ref{sec:weight-classes}, discuss Gaussian
nodes and weighted MZ inequalities in
Subsection~\ref{sec:gaussian-nodes}, and recall known
error estimates in Subsection~\ref{sec:known-estimates}.
Steklov averaging operators and their properties are
considered in Subsection~\ref{sec:steklov}, while
moduli of smoothness and $K$-functionals are discussed
in Subsection~\ref{sec:moduli-k}.
In Section~3, we present our main results,
dividing them into approximation by Lagrange
interpolation polynomials (Subsection~3.1)
and approximation by Kantorovich-type operators
(Subsection~\ref{secKant}).
Section~4 contains several auxiliary results
on the interval $[-1,1]$ (Subsection 4.1)
and on the torus $\mathbb T$ (Subsection 4.2),
as well as relations between various
$K$-functionals and moduli of smoothness
(Subsection 4.3).
 Finally, Section 5 is devoted to the proofs of the main results.

\smallskip

Throughout the paper, $c$, $C$, and $C_j$, $j=1,2,\dots$, denote positive constants which may vary from line to line and are independent of the essential parameters, typically $n$ and $f$. We write $A\asymp B$ if $C^{-1}A\le B\le C A$.

\section{Preliminary definitions and results}

Let $I := [-1,1]$. Consider a non-negative weight function $w: J \to \mathbb{R}$, where $J \subset I$, and let $1 \le p < \infty$. We define the weighted $L_p$-norm of a function $f$ on $J$ as
$$
\|f\|_{p, w}(J) := \|f\|_{L_{p,w}(J)} := \left( \int_{J} |f(x)|^p\, w(x)\, dx \right)^{1/p}.
$$
Accordingly, $L_{p,w}(J)$ denotes the space of all measurable functions $f: J \to \mathbb{R}$ for which $\|f\|_{p, w}(J) < \infty$. For simplicity, we drop the domain from the notation when $J = I$, writing $\|f\|_{p, w} = \|f\|_{p, w}(I)$.
By $L_p(I)$ and $\|\cdot\|_p$,
we denote the usual Lebesgue space
and its norm, respectively.
Moreover, $C_{\rm loc}(I)$ denotes the
class of functions on $I$ that are
continuous on $(-1,1)$.

We denote by ${\cal P}_n$ the set of all algebraic polynomials of degree at most $n$ on the interval $I$.
For $f\in L_{p,w}(I)$, we set
\[
E_n(f)_{p,w}:=
\inf_{P\in\mathcal P_n}\|f-P\|_{p,w}.
\]

Let $(p_n(u))_{n \in \Z_+}$ be the sequence of orthonormal polynomials with respect to a fixed weight $u$, where  $p_n\in {\cal P}_n$ and each polynomial is normalized to have a positive leading coefficient.
Denote by ${\cal L}_n(f, u) \in {\cal P}_{n-1}$ the Lagrange polynomial interpolating a function $f:I \to \R$ at the zeros $-1 < x_{1n}(u) < x_{2n}(u) < \dots < x_{nn}(u) < 1$ of $p_n(u)$. Specifically, we write:
\be \label{general-Lagrange}
{\cal L}_n(f, u)(x) := \sum_{k=1}^n f(x_k)\,\ell_{k, n}(x), \quad \quad \ell_{k,n}(x) :=
\frac{p_n(u,x)}
     {p_n'(u,x_k)(x-x_k)}.
\ee
Here and hereafter,
$$
x_k=x_{kn}=x_{kn}(u),\quad k=1,\dots,n.
$$
For convenience, we also denote
$$
x_0:=-1\quad\text{and}\quad x_{n+1}:=1.
$$

Since we study the Lagrange polynomials in the weighted $L_{p}$-setting, ${\cal L}_n(f, u)$ must be well-defined. Therefore, we treat $L_{p,w}(I)$ as a space of finite-valued, pointwise-defined functions rather than equivalence classes of elements defined almost everywhere.

\subsection{Classes of weight functions}\label{sec:weight-classes}

We use the class $J_p^*$ of weights introduced
by Ditzian and Totik~\cite[Definition~8.1.1]{DITO1987}.

\begin{definition}
Let $1\le p<\infty$. A positive weight $v$
on $(-1,1)$ belongs to $J_p^*$ if
\[
v(x)=v_-(\sqrt{1+x})v_+(\sqrt{1-x}),
\]
where
\[
v_\pm(t)=t^{\gamma_\pm}L_\pm(t),
\qquad \gamma_\pm>-2/p,
\]
and $L_\pm$ are positive functions bounded
above and below on every interval
$[\delta,\sqrt{2}]$, $\delta>0$, such that,
for every $\varepsilon>0$,
\[
t^\varepsilon L_\pm(t)
\quad\text{are nondecreasing, while}\quad
t^{-\varepsilon}L_\pm(t)
\quad\text{are nonincreasing}
\]
on $(0,\delta_\varepsilon)$ for some
$\delta_\varepsilon>0$.
\end{definition}

Recall that the weight $w$ satisfies the \emph{doubling condition} if there exists a constant $C>0$ such that
\begin{equation}\label{doub}
w(2J) \le C\, w(J),
\end{equation}
for all intervals $J \subset I$, where $2J$ denotes the interval concentric with $J$ and of twice its length.  As usual, for any measurable set $E \subset I$, we define
$w(E):= \int_E w(x)dx$.  In what follows, all integrals over a subset $E$ are understood as integrals over $E \cap I$.

For a weight $w$ on $I$, we define the associated even
$2\pi$-periodic weight by
\begin{equation*}
W(t) := w(\cos t)\, |\sin t|, \quad t\in\T:=(-\pi,\pi].
\end{equation*}

\begin{definition}
The weight $W$ satisfies Muckenhoupt's $A_p(\T)$ condition with $1<p<\infty$ if
\begin{equation}\label{ApI}
\sup_J
\left( \frac{1}{|J|} \int_J W(t)dt \right)
\left( \frac{1}{|J|} \int_J W(t)^{-\frac{1}{p-1}}dt \right)^{p-1}<\infty,
\end{equation}
where the supremum is taken over all possible intervals $J$ in $\T$ and $|J|$ is the length of $J$.
\end{definition}

If $w$ is a doubling weight on $I$,
then the associated periodic weight $W$
is doubling on $\T$.
The stronger condition $W\in A_p(\T)$
will be assumed explicitly whenever needed.

\subsection{Gaussian nodes and weighted Marcinkiewicz--Zygmund inequalities}\label{sec:gaussian-nodes}

Let $u$ be a positive integrable weight on
$I=[-1,1]$, and let
\[
X_n=\{x_{kn}(u):k=1,\ldots,n\},
\qquad x_{kn}(u)=\cos\theta_{kn},
\]
be the zeros of the orthogonal polynomial
of degree $n$ associated with $u$.
We set $\theta_{0n}=\pi$ and
$\theta_{n+1,n}=0$.

For a positive integrable weight $w$,
the algebraic Christoffel function is defined by
\begin{equation}\label{christ}
\lambda_n(w,x):=
\inf_{\substack{P\in\mathcal P_n\\ |P(x)|=1}}
\int_{-1}^{1}|P(t)|^2w(t)\,dt.
\end{equation}

Our main results are formulated in terms
of the following two conditions on the
interpolation nodes and weights.
The first condition concerns the angular
distribution of the Gaussian nodes:
\begin{equation}\tag{S}\label{eqshi2}
\frac{c_1}{n}\le
\theta_{kn}-\theta_{k+1,n}
\le\frac{c_2}{n},
\qquad k=0,\ldots,n,
\end{equation}
where $c_1,c_2>0$ are independent of $n$.
The second condition is the weighted
polynomial Marcinkiewicz--Zygmund inequality:
for $1<p<\infty$,
\begin{equation}\tag{MZ}\label{eqmzp}
\|P\|_{p,w}\asymp
\left(
\sum_{k=1}^{n}
\lambda_n(w,x_{kn}(u))
|P(x_{kn}(u))|^p
\right)^{1/p},
\qquad P\in\mathcal P_{n-1},
\end{equation}
with constants independent of $n$ and $P$.

Formulating our results under
\eqref{eqshi2} and~\eqref{eqmzp},
rather than restricting the weights
to a particular class, allows us to
treat quite general interpolation settings.
In particular, condition~\eqref{eqshi2}
holds for Gaussian nodes associated with
any doubling weight $u$;
see~\cite{MATO2010}.
Weighted polynomial MZ inequalities
have been studied for various weights
and node systems;
see~\cite{MARU1999,MATO2000,KOTI2026}.
The applicability of these results to
\eqref{eqmzp} depends on the particular
pair of weights $u$ and $w$.

To illustrate the assumptions~\eqref{eqshi2} and~\eqref{eqmzp},
we consider the Jacobi weights
\[
u(x)=(1-x)^\alpha(1+x)^\beta,
\qquad
w(x)=(1-x)^a(1+x)^b,
\]
where $\alpha,\beta,a,b>-1$, and
introduce the parameter condition
\begin{equation}\tag{J}\label{JacMZ}
\frac14-\frac1p
<\frac ap-\frac{\alpha}{2}
<\frac54-\frac1p,
\qquad
\frac14-\frac1p
<\frac bp-\frac{\beta}{2}
<\frac54-\frac1p.
\end{equation}
Equivalently, these conditions can be written as
\[
\frac{w^{1/p}}{\sqrt{u\varphi}}\in L_p(I),
\qquad
\frac{\sqrt{u\varphi}}{w^{1/p}}\in L_{p'}(I).
\]
It is known that under condition~\eqref{JacMZ}, the Gaussian nodes
associated with $u$ and the weight $w$
satisfy~\eqref{eqmzp}.
Indeed, the direct MZ estimate follows
from~\cite[Theorem~4.2]{MATO2000},
while the reverse estimate follows
from~\cite[Theorem~2.7,
(2.23)--(2.24)]{MARU1999}.

\subsection{Known error estimates}\label{sec:known-estimates}

Recall the definition of $r$th weighted
main-part modulus of smoothness of a function $f \in L_{p,w}(I)$  (see \cite{DITO1987}):
\begin{equation}\label{eqmainpartmod}
\Omega^r_{\varphi}(f,t)_{p,w} := \sup_{0 < h \le t} \| \Delta^r_{h\varphi} f \|_{p, w}(I_{rh}),
\end{equation}
where
$$
\varphi(x) := \sqrt{1 - x^2},
$$
\[
\Delta^r_{h\varphi} f(x) := \sum_{i=0}^{r} (-1)^i \binom{r}{i}
f\left(x + \frac{r h}{2}\varphi(x) - i h \varphi(x)\right), \quad r \ge 1,
\]
and
\[
I_{rh} := \left[-1 + 2h^2 r^2,\, 1 - 2h^2 r^2\right].
\]

Mastroianni and Russo~\cite[Theorem~3.1]{MARU1999}
established, under suitable assumptions on
the weights $u$ and $w$, the
following direct estimate for Lagrange interpolation
at Gaussian nodes:
\begin{equation}\label{eq:MR1999}
\|f-\mathcal L_n(f,u)\|_{p,w}
\le \frac{C}{n^{1/p}}
\int_0^{1/n}
\frac{\Omega_\varphi^\nu(f,t)_{p,w}}
{t^{1+1/p}}\,dt,
\end{equation}
where $1<p<\infty$, $\nu\in\mathbb N$, and
$f\in C_{\rm loc}(I)$ satisfies
\[
\int_0^1
\frac{\Omega_\varphi^\nu(f,t)_{p,w}}
{t^{1+1/p}}\,dt<\infty.
\]
For functions with derivatives of bounded variation,
they also established in~\cite{MARU2010}
the following estimate for Lagrange interpolation
with suitable Jacobi weights $u$ and $w$:
\begin{equation}\label{eq:MR2010}
\|f-\mathcal L_n(f,u)\|_{p,w}
\le \frac{C}{n^{r+1/p}}
\int_{-1}^1
\varphi(t)^{r+1/p}w(t)^{1/p}
\,|df^{(r)}(t)|,
\qquad n>r,
\end{equation}
where $1<p<\infty$, $r\in\mathbb N_0$, and
$f^{(r)}\in BV(I)$.

\vskip0.2cm

In this paper, we extend and refine the aforementioned error estimates by considering a broader class of functions and employing measures of smoothness based on the Steklov averaging operators. Moreover, we establish the corresponding inverse error estimates, which, together with the direct ones, yield a complete characterization of the rate of convergence in terms of Steklov-based measures of smoothness introduced below.

\subsection{Steklov averaging operators and their properties}\label{sec:steklov}

In this subsection, we introduce and study Steklov-type averaging operators adapted to the geometry
of the interval $[-1,1]$.

\begin{definition}\label{steklov}
Let $h>0$ and suppose that
$f(\cos\cdot)\in L_1(\T)$.
The symmetric Steklov averaging operator
is defined by
\[
\mathcal A_h f(x)
:=
\frac{1}{2h}
\int_{\arccos x-h}^{\arccos x+h}
f(\cos v)\,dv,
\qquad x\in I.
\]
\end{definition}

As usual, for $r\in\mathbb N$, we set
\begin{equation*}
(\I-\mathcal A_h)^r f
=
\sum_{k=0}^r
(-1)^k\binom{r}{k}\mathcal A_h^k f,
\end{equation*}
where $\mathcal A_h^k$ denotes the
$k$-fold iterate of $\mathcal A_h$,
with $\mathcal A_h^0=\I$.

For the Gaussian nodes $X_n$ associated with
the weight $u$ (see Subsection~\ref{sec:gaussian-nodes}),
we define the discrete weighted norm
\begin{equation*}
\|f\|_{\ell_{p,w}(X_n)}
:=
\bigg\|
\sum_{k=1}^n |f(x_k)|
\chi_{[x_k,x_{k+1}]}
\bigg\|_{p,w}.
\end{equation*}
It is easy to see that this norm can be explicitly written as
\begin{equation}\label{discrnormlp}
  \|f\|_{\ell_{p,w}(X_n)}=\(\sum_{k=1}^n \mu_k|f(x_k)|^p\)^{1/p},
\end{equation}
where the weights are defined by
\begin{equation*}
  \mu_k := \int_{x_{k}}^{x_{k+1}} w(x) \, dx, \qquad k=1,\ldots, n.
\end{equation*}
See also Lemma~\ref{lemour} for an equivalent formula for $\|f\|_{\ell_{p,w}(X_n)}$.

The main results of this paper involve
symmetric Steklov differences in both
continuous and discrete weighted norms.
Below, we establish several estimates
relating the corresponding quantities to each other
and to more conventional measures
of smoothness.

\begin{theorem}\label{thbestapp}
Let $1<p<\infty$, $r\in\mathbb N$,
$\gamma>0$, and let $u$ and $w$ be
positive integrable weights satisfying
conditions~\eqref{eqshi2} and~\eqref{eqmzp},
with $W\in A_p(\T)$.
For $f\in L_{p,w}(I)$ and $n\in\mathbb N$,
suppose that $P_n\in\mathcal P_{n-1}$
satisfies
$\|f-P_n\|_{p,w}\le 2E_{n-1}(f)_{p,w}$.
Then
\begin{equation}\label{eqomega0}
       \| ({\rm I} - {\cal A}_{\gamma/n})^r f\|_{\ell_{p,w}(X_n)}\le C\(\| ({\rm I} - {\cal A}_{\gamma/n})^r f\|_{p,w}+\|f-P_n\|_{\ell_{p,w}(X_n)}\),
    \end{equation}
where the constant $C$ is independent of $f$ and $n$.
\end{theorem}

The preceding theorem estimates the discrete
norm of the Steklov differences in terms of
their continuous norm and the polynomial
approximation error at the nodes.
The next result provides a related estimate
involving the weighted Ditzian--Totik
main-part modulus of smoothness.

\begin{theorem}\label{thOmega}
Let $1<p<\infty$, $\nu,r\in\mathbb N$,
and $\gamma>0$.
Suppose that $u$ and $w$ are positive
integrable weights such that
condition~\eqref{eqshi2} holds,
$w^{1/p}\in J_p^*$, and $W\in A_p(\T)$.
Then, for every
$f\in L_{p,w}(I)\cap C_{\rm loc}(I)$
satisfying
\[
\int_0^1
\frac{\Omega_\varphi^\nu(f,t)_{p,w}}
{t^{1+1/p}}\,dt<\infty,
\]
we have
\begin{equation}\label{intOmega}
\begin{split}
\|(\I-\A_{\gamma/n})^r f\|_{\ell_{p,w}(X_n)}
\le C\bigg(
&\|(\I-\A_{\gamma/n})^r f\|_{p,w}+n^{-1/p}\int_0^{1/n}
\frac{\Omega_\varphi^\nu(f,t)_{p,w}}
{t^{1+1/p}}\,dt
\bigg),
\end{split}
\end{equation}
where $C$ is independent of $f$ and $n$.
\end{theorem}

We now recall the notion of bounded
$p$-variation in the sense of Wiener.
For $p\ge 1$, let
\[
WBV_p(I):=
\{f: V_p(f,I)<\infty\},
\]
where
\[
V_p(f,I):=
\sup_D\sum_{i=1}^m
|f(t_i)-f(t_{i-1})|^p,
\]
and the supremum is taken over all
partitions $D=\{t_0,t_1,\ldots,t_m\}$
of $I$ (see, e.g.,~\cite{ABM}).

The corresponding $p$-variation function
is defined by
\[
v_p(f,x):=V_p(f,[-1,x]),
\qquad x\in I.
\]
We use $dv_p(f,x)$ to denote the
Lebesgue--Stieltjes measure associated
with the right-continuous representative
of $v_p(f,\cdot)$.

\begin{theorem}\label{th2-open-problem2}
Let $1<p<\infty$, and let
$w(x)=(1-x)^a(1+x)^b$, where $a,b>-1$.
Suppose that the Gaussian nodes associated
with $u$ satisfy~\eqref{eqshi2}.
Let $f\in C_{\rm loc}(I)\cap WBV_p(I)$
and assume that $\gamma\in(0,\pi/2)$
satisfies
\begin{equation}\label{gprob2}
\gamma<
\inf_{n\ge1}n\min_{0\le k\le n}
(\theta_{kn}-\theta_{k+1,n}).
\end{equation}
Then
\[
\|f-\mathcal A_{\gamma/n}f\|_{\ell_{p,w}(X_n)}
\le \frac{C}{n^{1/p}}
\left(
\int_{-1+c/n^2}^{1-c/n^2}
\varphi(x)w(x)\,dv_p(f,x)
\right)^{1/p},
\]
where $C,c>0$ are independent of $f$ and $n$.
\end{theorem}

\begin{remark}\label{remproblem2}
    Note that, proceeding similarly to the proof of Theorem~\ref{th2-open-problem2} (see Subsection~\ref{secth2-open-problem2}), we obtain an analogous result for functions in the classical $BV$-space. Namely, by applying the triangle inequality first and then using the relation $dv_1(x) = |df(x)|$ for the total variation function $v_1$, for all $f \in C_{\rm loc}(I) \cap BV(I)$ we have
$$
\|f-{\cal A}_{\gamma/n} f\|_{\ell_{p,w}(X_n)}\le \frac{C}{n^{1/p}}\, \int_{-1}^1 (\vp(x) w(x))^{1/p}\, |d f(x)|,\, \quad n \in \mathbb{N}.
$$
\end{remark}

\begin{theorem}\label{th2-open-problem2+}
Let $1<p<\infty$, $\gamma>0$, and let $w$
be a positive integrable doubling weight on $I$.
Then, for every
$f\in C_{\rm loc}(I)\cap WBV_p(I)$ and
$n\in\mathbb N$, we have
\begin{equation}\label{eqeq1}
\|f-\mathcal A_{\gamma/n}f\|_{p,w}
\le C
\left(
\int_{-1}^1
\int_{I\cap[x-\Delta_n(x),x+\Delta_n(x)]}
w(y)dy\,dv_p(f,x)
\right)^{1/p},
\end{equation}
where
\[
\Delta_n(x):=\frac{\varphi(x)}n+\frac1{n^2},
\]
and $C$ is independent of $f$ and $n$.
\end{theorem}

\begin{remark}
In general, the estimate~\eqref{eqeq1}
cannot be simplified by replacing the
local integral of $w$ with
$n^{-1}\varphi(x)w(x)$.
Indeed, the estimate
\begin{equation}\label{eqeq2}
\|f-\mathcal A_{\gamma/n}f\|_{p,w}
\le \frac{C}{n^{1/p}}
\left(
\int_{-1}^1
\varphi(x)w(x)\,dv_p(f,x)
\right)^{1/p}
\end{equation}
fails even for continuous functions.

To see this, take $w(x)=(1-x)^p$,
fix $n$ sufficiently large, and put
$h=\gamma/n<\pi/4$.
For $0<\varepsilon<h/8$, define
\[
f_\varepsilon(\cos\theta)
=\psi(\theta/\varepsilon),
\qquad 0\le\theta\le\pi,
\]
where
\[
\psi(t)=
\begin{cases}
1, & 0\le t\le 1/2,\\
2(1-t), & 1/2<t<1,\\
0, & t\ge1.
\end{cases}
\]
Then $f_\varepsilon$ is continuous and
monotone, and its $p$-variation measure
is supported in
$[\cos\varepsilon,\cos(\varepsilon/2)]$,
with total mass one.
Since
$w(x)\varphi(x)\le C\varepsilon^{2p+1}$
on this interval, we obtain
\[
\int_{-1}^1
\varphi(x)w(x)\,dv_p(f_\varepsilon,x)
\le C\varepsilon^{2p+1}.
\]
On the other hand, for
$\theta\in[h/2,3h/4]$, we have
$f_\varepsilon(\cos\theta)=0$ and
\[
\mathcal A_hf_\varepsilon(\cos\theta)
=\frac{\varepsilon}{h}
\int_0^1\psi(t)\,dt
=\frac{3\varepsilon}{4h}.
\]
Hence
\[
\begin{aligned}
\|f_\varepsilon-\mathcal A_hf_\varepsilon\|_{p,w}^p
&\ge
\left(\frac{3\varepsilon}{4h}\right)^p
\int_{h/2}^{3h/4}
(1-\cos\theta)^p\sin\theta\,d\theta\\
&\ge c_h\varepsilon^p,
\end{aligned}
\]
where $c_h>0$ is independent of
$\varepsilon$.
Thus the left-hand side of~\eqref{eqeq2}
is bounded below by $c_h^{1/p}\varepsilon$,
whereas its right-hand side is
$O(\varepsilon^{2+1/p})$.
Letting $\varepsilon\to0$ gives
a contradiction.
\end{remark}

\subsection{Moduli of smoothness and $K$-functionals}\label{sec:moduli-k}

In this subsection, we introduce several
quantities used to characterize
Steklov differences in weighted
$L_p$-spaces. For simplicity, we restrict ourselves to the case of the Jacobi weights.
Let
\[
w^{(a,b)}(x):=(1-x)^a(1+x)^b,
\qquad a,b>-1,
\]
denote the Jacobi weight. We write
$$
 L_{p,(a,b)}:=L_{p,w^{(a,b)}}(I) \quad\text{and}\quad \|\cdot\|_{p,(a,b)}:=\|\cdot\|_{p,w^{(a,b)}}.
$$

Recall (see~\cite[(2.5.3)]{MAMI2008}) that the (global) weighted modulus of smoothness of order $r \ge 1$ is given by
\begin{equation*}
  \begin{split}
     \omega_{\varphi}^r(f,t)_{p,(a,b)}
:= \Omega_{\varphi}^r(f,t)_{p,(a,b)}
+ \inf_{P \in \mathcal{P}_{r-1}} &\|f-P\|_{L_{p,(a,b)}\([-1,\,-1+4{r^2} t^2]\)}\\
&+ \inf_{P \in \mathcal{P}_{r-1}} \|f-P\|_{L_{p,(a,b)}\([1-4{r^2} t^2,\,1]\)}.
   \end{split}
\end{equation*}
To establish the connection with the corresponding $K$-functional, we consider the pair of spaces
$({L_{p,(a,b)}},\, W_{p,{(a,b)}}^r)$, where the Sobolev-type space $W_{p,{(a,b)}}^r$ consists of functions $g\in L_{p, (a,b)}$ such that $g^{(r-1)}\in AC_{\rm loc}(I)$ and $\varphi^rg^{(r)}\in L_{p,(a,b)}$.
For $r \in \mathbb{N}$, the $K$-functional $K_{\varphi}^r(f,t)_{p,(a,b)}$ is then defined as
$$
K_{\varphi}^r(f,t)_{p,(a,b)}
:= \inf\left\{
\|f-g\|_{p,(a,b)} + t^r \|\varphi^r g^{(r)}\|_{p,(a,b)}
: g \in W_{p,{(a,b)}}^r
\right\}.
$$

Certain relations between these moduli of smoothness and $K$-functionals are discussed in detail in Subsection~\ref{secMK}.

To characterize  Steklov differences
of higher orders, we use the differential
operators and the corresponding
$K$-functionals associated with
Fourier--Chebyshev expansions.

Let $C_k$, $k=0,1,\dots$, denote the orthonormal Chebyshev polynomials of the first kind on $I$ with respect to the weight $w^{(-1/2,-1/2)}$, defined explicitly as
$$
C_0(x)=\frac{1}{\sqrt{\pi}},\qquad
C_k(x)=\sqrt{\frac{2}{\pi}}\cos(k\arccos x),\quad k=1,2,\dots.
$$
These polynomials are the eigenfunctions of the differential operator
$$
\mathcal{D}_2
:= -\left(\varphi(x) \frac{d}{dx}\right)^2,
$$
satisfying $\mathcal{D}_2 C_k = k^2 C_k$ for all $k = 0, 1, 2, \dots$.

For a function $f \in L_{p,{(-1/2, -1/2)}}$ with $1 \le p \le \infty$, its Fourier--Chebyshev expansion is given by
$$
f(x) \sim \sum_{k=0}^{\infty} \hat{f}_k C_k(x),\qquad \text{where} \quad
\hat{f}_k
:= \int_{-1}^{1} f(x) C_k(x) \,\frac{dx}{\sqrt{1-x^2}},
\qquad k = 0, 1, 2, \dots .
$$
For any $\sigma > 0$, if there exists a function $g \in L_{1,{(-1/2, -1/2)}}$ whose Fourier--Chebyshev expansion is of the form
$$
g \sim \sum_{k=1}^{\infty} k^{\sigma} \hat{f}_k C_k,
$$
then we write $g = \mathcal{D}_{\sigma} f$, and we call $\mathcal{D}_{\sigma} f$ the fractional derivative of order $\sigma$ of the function $f$.

We now assume that $1 <p < \infty$ and
\begin{equation}\label{condGT}
-1<a<\frac p2-1, \qquad -1<b<\frac p2-1.
\end{equation}
Under these conditions, since the embedding $L_{p,{(a,b)}} \subset L_{1,{(-1/2,-1/2)}}$ holds, the Fourier--Che\-byshev expansion is well-defined for any $f \in L_{p,{(a,b)}}$. The $K$-functional corresponding to the fractional differential operator $\mathcal{D}_r$ for $r > 0$ is defined as
$$
{\rm K}_r(f, t)_{p,(a,b)}
:= \inf\left\{
\|f - g\|_{p,(a,b)} + t^r \|\mathcal{D}_r g\|_{p,(a,b)} :
g \in {\rm W}_{p,(a,b)}^{r}
\right\},
$$
where the underlying space is given by
$$
{\rm W}_{p,(a,b)}^{r}
:= \{\, g : g,\, \mathcal{D}_r g \in L_{p,{(a,b)}} \,\}.
$$

\begin{theorem}\label{thArK}
Let $1<p<\infty$, $-1<a,b<p/2-1$,
$0<h<1$, and $r\in\N$. Then, for every
$f\in L_{p,(a,b)}$, we have
\begin{equation}\label{eqArK}
\|(\I-\A_h)^r f\|_{p,(a,b)}
\asymp {\rm K}_{2r}(f,h)_{p,(a,b)}.
\end{equation}
\end{theorem}

\section{Main results}

\subsection{Approximation by Lagrange interpolation polynomials}

We start with the direct error estimate.

\begin{theorem}\label{thsp}
Let $1<p<\infty$, $r\in\mathbb N$,
and $\gamma>0$.
Suppose that $u$ and $w$ are positive
integrable weights satisfying
conditions~\eqref{eqshi2} and~\eqref{eqmzp},
with $W\in A_p(\T)$.
Then, for every $f\in L_{p,w}(I)$,
we have
\begin{equation}\label{stp1}
\begin{split}
\|f-\mathcal L_n(f,u)\|_{p,w}
\le C\bigl(
&\|(\I-\mathcal A_{\gamma/n})^r f\|_{p,w}+\|(\I-\mathcal A_{\gamma/n})^r f
\|_{\ell_{p,w}(X_n)}
\bigr),
\end{split}
\end{equation}
where $C$ is independent of $f$ and $n$.
\end{theorem}

The next theorem provides the corresponding inverse error estimates.

\begin{theorem}\label{prinvp}
Under the assumptions of Theorem~\ref{thsp},
we have
\begin{equation}\label{st1invp}
\begin{split}
C_1\|(\I-\A_{\g/n})^r f\|_{\ell_{p,w}(X_n)}
-C_2\|(\I-\A_{\g/n})^r f\|_{p,w}
\le \|f-\mathcal L_n(f,u)\|_{p,w}
\end{split}
\end{equation}
and
\begin{equation}\label{inv0p}
\begin{split}
\|(\I-\A_{\g/n})^r f\|_{p,w}
&+\|(\I-\A_{\g/n})^r f\|_{\ell_{p,w}(X_n)}
\\
&\le C\left(
\|f-\mathcal L_n(f,u)\|_{p,w}
+\frac1{n^{2r}}
\sum_{k=1}^n k^{2r-1}
\|f-\mathcal L_k(f,u)\|_{p,w}
\right),
\end{split}
\end{equation}
where the constants are independent of
$f$ and $n$.
\end{theorem}

The application of Theorems~\ref{thsp} and~\ref{prinvp} yields the following corollaries.

\begin{corollary}\label{cor2}
Under the assumptions of Theorem~\ref{thsp},
the following properties are equivalent
for every $\lambda\in(0,2r)$:
\begin{itemize}
\item[(i)]
$\|f-\mathcal L_n(f,u)\|_{p,w}
=\mathcal O(n^{-\lambda})$;

\item[(ii)]
$\|(\I-\A_{\g/n})^r f\|_{p,w}
+\|(\I-\A_{\g/n})^r f\|_{\ell_{p,w}(X_n)}
=\mathcal O(n^{-\lambda})$.
\end{itemize}
\end{corollary}

Taking also into account that
$\|(\I-\A_{\g/n})^r f\|_{p,w}\to0$
as $n\to\infty$ for every
$f\in L_{p,w}(I)$, we obtain
the following convergence criterion.

\begin{corollary}\label{cor1}
%Let $f\in C_{\rm loc}(I)$.
Under the assumptions of Theorem~\ref{thsp}, we have
  $$
    \|f-{\cal L}_n(f,u)\|_{p,w}\to 0\quad\text{as } n\to \infty
  $$
if and only if
\begin{equation*}%\label{cor1.1}
  \|(\I-\A_{\g/n})^r f\|_{\ell_{p,w}(X_n)}\to 0\quad\text{as } n\to \infty.
\end{equation*}
\end{corollary}

Combining Theorems~\ref{thsp},
\ref{th2-open-problem2}, and
\ref{th2-open-problem2+}, we obtain the
following error estimate for functions
of bounded $p$-variation.
For the Jacobi weight
$w(x)=(1-x)^a(1+x)^b$, $a,b>-1$, put
\[
w_n(x):=
\left(\sqrt{1-x}+\frac1n\right)^{2a}
\left(\sqrt{1+x}+\frac1n\right)^{2b}.
\]

\begin{corollary}
Let $1<p<\infty$, and let
$w(x)=(1-x)^a(1+x)^b$, where $a,b>-1$
and $W\in A_p(\mathbb T)$.
Assume that the Gaussian nodes associated
with $u$ and the weight $w$ satisfy
conditions~\eqref{eqshi2} and~\eqref{eqmzp}.
Let $0<\gamma<\pi/2$ satisfy~\eqref{gprob2}.
Then, for every
$f\in C_{\rm loc}(I)\cap WBV_p(I)$
and $n\in\mathbb N$, we have
\begin{equation*}
\begin{split}
\|f-\mathcal L_n(f,u)\|_{p,w}
&\le C\|f-\mathcal A_{\gamma/n}f\|_{p,w}+\frac{C}{n^{1/p}}
\left(
\int_{-1+c/n^2}^{1-c/n^2}
\varphi(x)w(x)\,dv_p(f,x)
\right)^{1/p}\\
&\le\frac{C}{n^{1/p}}
\left(
\int_{-1}^1
\left(\varphi(x)+\frac1n\right)
w_n(x)\,dv_p(f,x)
\right)^{1/p},
\end{split}
\end{equation*}
where
the constants $C,c>0$ are independent
of $f$ and $n$.
\end{corollary}

The second inequality follows from
\begin{equation}\label{eq:jacobi-local-weight}
\int_{I\cap[x-\Delta_n(x),x+\Delta_n(x)]}
w(y)\,dy
\asymp \Delta_n(x)w_n(x),
\end{equation}
which is a consequence of~\eqref{eqmato1}
and~\cite[Lemma~5, p.~108]{NEVAI1979},
together with
\[
w(x)\le Cw_n(x),
\qquad
x\in[-1+cn^{-2},1-cn^{-2}].
\]
The latter estimate is a special case
of~\cite[(8.1.5)]{DITO1987}.

We next obtain an estimate of the form~\eqref{eq:MR1999}
under the assumptions~\eqref{eqshi2} and~\eqref{eqmzp}.

\begin{corollary}\label{corthsp}
Let $1<p<\infty$, $\nu\in\mathbb N$, and suppose
that $w^{1/p}\in J_p^*$.
Assume that conditions~\eqref{eqshi2}
and~\eqref{eqmzp} hold.
Then, for all sufficiently large $n$ and every
$f\in L_{p,w}(I)\cap C_{\rm loc}(I)$ satisfying
\[
\int_0^1
\frac{\Omega_\varphi^\nu(f,t)_{p,w}}
{t^{1+1/p}}\,dt<\infty,
\]
we have
\begin{equation}\label{eq:lagrange-jackson}
\|f-\mathcal L_n(f,u)\|_{p,w}
\le Cn^{-1/p}
\int_0^{1/n}
\frac{\Omega_\varphi^\nu(f,t)_{p,w}}
{t^{1+1/p}}\,dt.
\end{equation}
\end{corollary}

\subsubsection{Some corollaries in the case of Jacobi weights}
For the usual Jacobi weight
$$
w^{(a,b)}(x)=(1-x)^a(1+x)^b,\quad a,b>-1,
$$
the above results can be formulated using more conventional measures of smoothness.

\begin{corollary}\label{cor_prinvp}
Let $1<p<\infty$, $r\in\N$, and
$\gamma>0$.
Assume that $\alpha,\beta>-1$,
$-1<a,b<p/2-1$, and condition~\eqref{JacMZ}
holds.
Then, for every $f\in L_{p,(a,b)}$, we have
\[
\begin{aligned}
\|f-\mathcal L_n(f,w^{(\alpha,\beta)})\|_{p,(a,b)}
\le C\bigl(
&{\rm K}_{2r}(f,\gamma/n)_{p,(a,b)}+\|(\I-\A_{\gamma/n})^r f
\|_{\ell_{p,(a,b)}(X_n)}
\bigr)
\end{aligned}
\]
and
\[
\begin{aligned}
&{\rm K}_{2r}(f,\gamma/n)_{p,(a,b)}
+\|(\I-\A_{\gamma/n})^r f
\|_{\ell_{p,(a,b)}(X_n)}\\
&\quad\le C\left(
\|f-\mathcal L_n(f,w^{(\alpha,\beta)})\|_{p,(a,b)}
+\frac1{n^{2r}}\sum_{k=1}^n
k^{2r-1}
\|f-\mathcal L_k(f,w^{(\alpha,\beta)})\|_{p,(a,b)}
\right).
\end{aligned}
\]
The constants are independent of $f$ and $n$,
and the estimates hold for all sufficiently
large $n$, depending on $\gamma$.
\end{corollary}

\begin{proof}
Since every Jacobi weight is doubling,
condition~\eqref{eqshi2} follows
from~\cite{MATO2010}.
Condition~\eqref{JacMZ} implies
\eqref{eqmzp}, as explained in Subsection~\ref{sec:gaussian-nodes}.
Furthermore, $-1<a,b<p/2-1$ implies
$W\in A_p(\T)$.
Hence, the assertion follows from
Theorems~\ref{thsp}, \ref{prinvp},
and~\ref{thArK}.
\end{proof}

Combining Corollary~\ref{cor_prinvp}
with Lemmas~\ref{lemGT}
and~\ref{lemDTKM}, we obtain
the following equivalences.

\begin{corollary}\label{cor2J}
Under the assumptions of
Corollary~\ref{cor_prinvp}, the following
properties are equivalent for every
$\lambda\in(0,2r)$:
\begin{itemize}
\item[(i)]
$\|f-\mathcal L_n(f,w^{(\alpha,\beta)})\|_{p,(a,b)}
=\mathcal O(n^{-\lambda})$;

\item[(ii)]
${\rm K}_{2r}(f,\gamma/n)_{p,(a,b)}
+\|(\I-\A_{\gamma/n})^r f
\|_{\ell_{p,(a,b)}(X_n)}
=\mathcal O(n^{-\lambda})$.
\end{itemize}
If, additionally, $a,b\ge0$, then these
properties are equivalent to
\begin{itemize}
\item[(iii)]
$\omega_{\varphi}^{2r}(f,\gamma/n)_{p,(a,b)}
+\|(\I-\A_{\gamma/n})^r f
\|_{\ell_{p,(a,b)}(X_n)}
=\mathcal O(n^{-\lambda})$;

\item[(iv)]
$K_{\varphi}^{2r}(f,\gamma/n)_{p,(a,b)}
+\|(\I-\A_{\gamma/n})^r f
\|_{\ell_{p,(a,b)}(X_n)}
=\mathcal O(n^{-\lambda})$.
\end{itemize}
\end{corollary}

\subsection{Approximation by Kantorovich-type operators}\label{secKant}

The main aim of this section is to study Kantorovich-type sampling operators constructed using the Lagrange polynomials ${\cal L}_n(f,u)$.

\begin{definition}
Let $r,n\in\N$ and $\gamma>0$.
 For any $f$ such that
$f(\cos\cdot)\in L_1(\T)$, we define
$$
\K_{n,\g}^r(f,u)(x) := \sum_{k=1}^n \sum_{\nu=1}^r (-1)^{\nu+1}\binom{r}{\nu}\A_{\g/n}^\nu f(x_k)\, \ell_{k,n}(x).
$$
\end{definition}

\begin{theorem}\label{thK1}
    Under the conditions of Theorem~\ref{thsp},
for all $f\in L_{p,w}(I)$, we have
    \begin{equation}\label{eqK1}
        \|f-\K_{n,\g}^r (f,u)\|_{p,w}\asymp \|(\I-\A_{\g/n})^r f\|_{p,w}.
    \end{equation}
\end{theorem}

For Jacobi weights, Theorem~\ref{thK1}
together with Theorem~\ref{thArK}
and Lemmas~\ref{lemDTKM} and~\ref{lemGT}
yields the following result.

\begin{corollary}\label{prinvpK}
Let $1<p<\infty$, $r\in\N$, and
$\gamma>0$.
Assume that $\alpha,\beta>-1$,
$-1<a,b<p/2-1$, and condition~\eqref{JacMZ}
holds.
Then for any  $f\in L_{{p},{(a,b)}}$, we have
\begin{equation*}%\label{stp1}
  \|f-\K_{n,\g}^r(f,w^{(\a,\b)})\|_{p,(a,b)}\asymp {\rm K}_{2r}(f, \g/n)_{p,(a,b)}.
\end{equation*}
If, additionally, $a,b\ge 0$, then
\begin{equation*}%\label{stp1}
  C^{-1}\omega_{\varphi}^{2r}(f,\g/n)_{p,(a,b)}\le \|f-\K_{n,\g}^r(f,w^{(\a,\b)})\|_{p,(a,b)}\le C\(\omega_{\varphi}^{2r}(f,\g/n)_{p,(a,b)}+n^{-2r}\|f\|_{p,(a,b)}\).
\end{equation*}
The constants are independent of $f$ and $n$,
and the estimates hold for all sufficiently
large $n$, depending on $\gamma$.
\end{corollary}

\begin{proof}
The assumptions imply (S), (MZ), and
$W\in A_p(\T)$, as in
Corollary~\ref{cor_prinvp}.
The assertion follows from
Theorem~\ref{thK1}, Theorem~\ref{thArK},
and Lemmas~\ref{lemDTKM} and~\ref{lemGT}.
\end{proof}

As usual, Corollary~\ref{prinvpK} implies the following equivalence results.

\begin{corollary}\label{prinvpK3}
Under the conditions of Corollary~\ref{prinvpK}, the following properties are equivalent for all $\lambda\in (0,2r)$:
  \begin{itemize}
    \item[(i)]  $\|f-\K_{n,\g}^r(f,w^{(\a,\b)})\|_{p,(a,b)}=\mathcal{O}(n^{-\lambda})$,
    \item[(ii)] $\|(\I-\A_{\g/n})^r f\|_{p,(a,b)}=\mathcal{O}(n^{-\lambda})$.
  \end{itemize}
 If, additionally, $a,b\ge 0$, the properties (i) and (ii)  are equivalent to
  \begin{itemize}
    \item[(iii)]  $\omega_{\varphi}^{2r}(f,\g/n)_{p,(a,b)}=\mathcal{O}(n^{-\lambda})$,
    \item[(iv)] $K_{\varphi}^{2r}(f,\g/n)_{p,(a,b)}=\mathcal{O}(n^{-\lambda})$.
  \end{itemize}
\end{corollary}

\section{Auxiliary results}

\subsection{Auxiliary results on $[-1,1]$} \label{sub-auxiliary1}

\begin{lemma}\label{lemchri} (See~\cite[(7.14)]{MATO2000}.)
Let $w$ satisfy the doubling condition.
Then, for every $n\in \N$ and $x\in I$, we have
\begin{equation}\label{eqmato1}
\lambda_n(w,x)\asymp\, \int_{x - \Delta_n(x)}^{x + \Delta_n(x)} w(t)\,dt \,.
\end{equation}
\end{lemma}

Recall the notation
$$
x_k=x_{kn}=x_{kn}(u),\quad k=1,\dots,n,\quad\text{and}\quad
x_0=x_{0n}=-1,\quad x_{n+1}=x_{n+1,n}=1,
$$
and
$$
\theta_{k}=\theta_{kn}=\arccos x_{kn}(u),\quad k = 0, 1, \ldots, n+1.
$$

\begin{lemma}\label{lemshi1}
Let $u$ be a doubling weight on $I$.
Then condition (S) holds. Moreover,
\begin{equation}\label{eqshi3}
c_3\Delta_n(x_{kn})
\le x_{k+1,n}-x_{kn}
\le c_4\Delta_n(x_{kn}),
\qquad k=0,\ldots,n,
\end{equation}
where $c_3,c_4>0$ are independent of $n$.
\end{lemma}

\begin{proof}
Condition~\eqref{eqshi2} follows
from~\cite{MATO2010}, while~\eqref{eqshi3}
follows from~\eqref{eqshi2}
and~\cite[Lemma~2.6]{SHI2005}.
\end{proof}

\begin{lemma}\label{lemour}
Let $w$ be a doubling weight and suppose
that condition~\eqref{eqshi2} holds. Then
\begin{equation}\label{eqour1}
  \lambda_n\!\left(w,\, x_{kn}\right)\asymp \mu_k
=\int_{x_{kn}}^{x_{k+1,n}}w(x)dx,\qquad k = 1, \ldots, n.
\end{equation}
In particular, for every  $f$ with finite
values at the interpolation nodes, we have
$$
\|f\|_{\ell_{p,w}(X_n)}\asymp \(\sum_{k=1}^{n}
|f(x_{kn})|^p\lambda_n(w,\, x_{kn})\)^{1/p}.
$$
\end{lemma}

\begin{proof}
Since the angular separation condition~\eqref{eqshi2}
implies the spacing estimates~\eqref{eqshi3},
we have
\begin{equation}\label{eqour3}
  \int_{x_{j}}^{x_{j+1}}w(x)dx
  \le
  \int_{x_j - c_4\Delta_n(x_j)}^{x_j + c_4\Delta_n(x_j)} w(x)dx,\quad j=1,\dots,n.
\end{equation}
Let $l\in \mathbb{N}$ be such that $c_4\le 2^l$. Then, using the following inclusion
$$
[x_j-c_4\Delta_n(x_j), x_j+c_4 \Delta_n(x_j)]\subseteq [x_j-2^l\Delta_n(x_j), x_j+2^l\Delta_n(x_j)]
$$
and applying $l$ times the doubling condition~\eqref{doub}, we get
\begin{equation*}
  \int_{x_j - c_4\Delta_n(x_j)}^{x_j + c_4\Delta_n(x_j)} w(x)dx
  \le
  C_l\int_{x_j - \Delta_n(x_j)}^{x_j + \Delta_n(x_j)} w(x)dx.
\end{equation*}
This together with~\eqref{eqour3} and Lemma~\ref{lemchri} implies the lower estimate in~\eqref{eqour1}.

To obtain the upper inequality, we set $m\ge 1+\frac{2}{c_3}$ and $d=x_{j+1}-x_j$. Then, using~\eqref{eqshi3}, we obtain
$$
\frac{x_{j}+x_{j+1}}{2}-\frac{md}{2}\le
\frac{x_{j}+x_{j+1}}{2}-\(1+\frac{2}{c_3}\)\frac{d}{2}=x_j-\frac{d}{c_3}
\le x_j-\Delta_n(x_j)
$$
and
$$
x_j+\Delta_n(x_j)\le x_j+\frac{d}{c_3}\le x_j+\(1+\frac{2}{c_3}\)\frac{d}{2}=\frac{x_{j}+x_{j+1}}{2}+\frac{md}{2}.
$$
Therefore, $[x_j - \Delta_n(x_j), x_j + \Delta_n(x_j)]\subset m [x_{j},x_{j+1}]$, which, together with repeated applications of the doubling property, yields
\begin{equation}\label{eqour5}
 \int_{x_j - \Delta_n(x_j)}^{x_j + \Delta_n(x_j)} w(x)dx
\le
 \int_{m[x_{j},x_{j+1}]} w(x)dx \le C_m \int_{x_j}^{x_{j+1}} w(x)dx.
\end{equation}
Finally, applying Lemma~\ref{lemchri}, we finish the proof.
\end{proof}

\begin{lemma}\label{lemmzf}
Let $1<p<\infty$, and let $w$ be
a doubling weight. Suppose that
conditions~\eqref{eqshi2} and~\eqref{eqmzp}
hold.
Then, for every $f$ with finite
values at the interpolation nodes, we have
\begin{equation}\label{asp2}
\|{\cal L}_n(f,u)\|_{p,w}
\asymp \|f\|_{\ell_{p,w}(X_n)}.
\end{equation}
\end{lemma}

\begin{proof}
Taking into account that ${\cal L}_n f\in \mathcal{P}_{n-1}$ and ${\cal L}_n f(x_k)=f(x_k)$, $k=1,\dots,n$, we obtain from Lemma~\ref{lemour}
and inequality~\eqref{eqmzp} that
\begin{equation*}%\label{asp2}
  \|f\|_{\ell_{p,w}(X_n)}=\|{\cal L}_n f\|_{\ell_{p,w}(X_n)}\asymp \(\sum_{k=1}^{n}
|{\cal L}_n f(x_{k})|^p\lambda_n(w,\, x_k)\)^{1/p} \asymp \|{\cal L}_n f\|_{p,w}.
\end{equation*}
\end{proof}

\begin{lemma}\label{lem:high-degree-sampling}
Let $1<p<\infty$ and
$w^{1/p}\in J_p^*$.
Suppose that condition~\eqref{eqshi2} holds.
Then, for all integers
$m\ge n\ge1$ and $P\in\mathcal P_m$,
\begin{equation}\label{eq:high-degree-sampling}
\|P\|_{\ell_{p,w}(X_n)}
\le C\left(\frac mn\right)^{1/p}
\|P\|_{p,w},
\end{equation}
where $C$ is independent of $P,m,n$.
\end{lemma}

\begin{proof}
Set $T(t)=P(\cos t)$ and
$W(t)=w(\cos t)|\sin t|$.
The defining conditions of $J_p^*$ imply
that $W$ is a doubling weight.
For $m\ge1$, put
\[
W_m(t)=m\int_{t-1/m}^{t+1/m}W(s)\,ds.
\]

By~\eqref{eqshi2}, the interval
$[\theta_{k+1},\theta_k]$ has length
comparable to $1/n$.
Hence, the local comparison
\cite[p.~93, (8.1.8)]{DITO1987},
applied to $w^{1/p}$, yields, for some
fixed $\eta>0$,
\[
W(t)\asymp W(\theta_k)
\quad\text{whenever}\quad
|t-\theta_k|\le\eta/m,\qquad m\ge n.
\]
Since $W$ is doubling, the above relation implies
\[
W_m(\theta_k)
\asymp m\int_{\theta_k-\eta/m}^{\theta_k+\eta/m}
W(t)\,dt
\asymp W(\theta_k),\qquad m\ge n.
\]
Similarly, by the doubling property of $W$
and~\eqref{eqshi2}, we obtain
\[
\begin{aligned}
\mu_k
&=\int_{\theta_{k+1}}^{\theta_k}W(t)\,dt\asymp
\int_{\theta_k-1/n}^{\theta_k+1/n}W(t)\,dt
=\frac1n W_n(\theta_k).
\end{aligned}
\]
Combining the preceding estimates, we get
\[
W_m(\theta_k)\asymp W_n(\theta_k)
\asymp n\mu_k,\qquad m\ge n.
\]

For $n\ge2$, applying the weighted large-sieve inequality
\cite[Theorem~2.3]{JOU2017}
to $|T|^2$ with exponent $p/2$
and degree at most $2m$, we obtain
\[
\sum_{k=1}^n
|T(\theta_k)|^p W_{2m}(\theta_k)
\le C(m+n)\|T\|_{p,W}^p
\le Cm\|T\|_{p,W}^p.
\]
For $n=1$, the same estimate follows by
applying the large-sieve inequality to
the set $\{\theta_1,\theta_1+\pi\}$.
Since $W$ is doubling,
$W_{2m}(t)\asymp W_m(t)$ uniformly in $t$.
Consequently,
\[
n\sum_{k=1}^n\mu_k|P(x_k)|^p
\le Cm\|P\|_{p,w}^p,
\]
which proves the assertion.
\end{proof}

\subsection{Auxiliary results on $\T$}

For every measurable function $f:I\to\mathbb C$,
we define its even $2\pi$-periodic lifting by
$F(t):=f(\cos t)$ and set
\[
\|F\|_{p,W}^*
:=
\left(\frac12\int_{-\pi}^{\pi}
|F(t)|^pW(t)\,dt\right)^{1/p}.
\]
Then
\begin{equation}\label{eqFf}
\|F\|_{p,W}^*=\|f\|_{p,w}.
\end{equation}
Similarly, for
$\Theta_n:=(\theta_k)_{k=1}^{n+1}$,
we define
\[
\|F\|_{\ell_{p,W}(\Theta_n)}
:=
\left(
\sum_{k=1}^n |F(\theta_k)|^p
\int_{\theta_{k+1}}^{\theta_k}W(t)\,dt
\right)^{1/p},
\]
so that
\begin{equation}\label{eqdFf}
\|F\|_{\ell_{p,W}(\Theta_n)}
=\|f\|_{\ell_{p,w}(X_n)}.
\end{equation}

For $F\in L_1^{\rm loc}(\T)$ and $h>0$,
the symmetric Steklov averaging operator
is defined by
\[
A_hF(t):=\frac1{2h}
\int_{t-h}^{t+h}F(v)\,dv.
\]
By the definition of $\mathcal A_h$,
\begin{equation}\label{eqFf+}
\mathcal A_hf(\cos t)=A_hF(t).
\end{equation}

Finally, we define the interpolation operator
$G_n$ by
\begin{equation}\label{eqFf+++}
G_nF(t)
:=
\sum_{k=1}^n
F(\theta_k)\ell_{k,n}(\cos t)
=\mathcal L_n(f,u)(\cos t).
\end{equation}

\begin{lemma}\label{lemmzF}
Let $1<p<\infty$ and $W\in A_p(\T)$.
Suppose that conditions~\eqref{eqshi2}
and~\eqref{eqmzp} hold.
Then
\begin{equation}\label{ast2}
\|G_nF\|_{p,W}^*
\asymp
\|F\|_{\ell_{p,W}(\Theta_n)}
\end{equation}
for every even function $F$ with finite
values at the interpolation nodes.
\end{lemma}

\begin{proof}
By~\eqref{eqFf} and~\eqref{eqFf+++},
\[
\|G_nF\|_{p,W}^*
=
\|\mathcal L_n(f,u)\|_{p,w}.
\]
Since $W\in A_p(\T)$, the weight $w$
is doubling on $I$. Thus, by the
interpolation property, the polynomial
inequality~\eqref{eqmzp},
and Lemma~\ref{lemour}, we obtain
\[
\|\mathcal L_n(f,u)\|_{p,w}
\asymp
\left(\sum_{k=1}^n
\mu_k|f(x_k)|^p\right)^{1/p}.
\]
The assertion now follows from~\eqref{eqdFf}.
\end{proof}

\begin{lemma}\label{lemsF}
(See~\cite[\S~5.2.1]{ST1993}.)
Let $1<p<\infty$ and $W\in A_p(\T)$.
Then, for every $F\in L_{p,W}(\T)$
and $h>0$, we have
\[
\|A_hF\|_{p,W}^*
\le C\|F\|_{p,W}^*.
\]
\end{lemma}

\begin{lemma}\label{lemsFd}
Let $1<p<\infty$, $W\in A_p(\T)$,
and suppose that condition~\eqref{eqshi2} holds.
If $h\asymp n^{-1}$, then
\[
\|A_hF\|_{\ell_{p,W}(\Theta_n)}
\le C\|F\|_{p,W}^*
\]
for every $F\in L_{p,W}(\T)$.
\end{lemma}

\begin{proof}
Let
$J_k:=[\theta_{k+1},\theta_k]$,
$k=1,\ldots,n$,
and put
$M:=\max_{1\le k\le n}
(\theta_k-\theta_{k+1}).$
By~\eqref{eqshi2}, we have $M\le Cn^{-1}$.
For every $t\in J_k$,
\[
[\theta_k-h,\theta_k+h]
\subset[t-(M+h),t+(M+h)].
\]
Consequently,
\[
|A_hF(\theta_k)|
\le\frac{M+h}{h}A_{M+h}|F|(t).
\]
Since $h\asymp n^{-1}$, the factor
$(M+h)/h$ is uniformly bounded.
It follows that
\[
\begin{aligned}
\|A_hF\|_{\ell_{p,W}(\Theta_n)}^p
&=\sum_{k=1}^n
|A_hF(\theta_k)|^p
\int_{J_k}W(t)\,dt\\
&\le C\sum_{k=1}^n
\int_{J_k}
|A_{M+h}|F|(t)|^pW(t)\,dt\\
&\le C\int_0^\pi
|A_{M+h}|F|(t)|^pW(t)\,dt\\
&\le C\|F\|_{p,W}^{*p},
\end{aligned}
\]
where the last inequality follows from
Lemma~\ref{lemsF}.
\end{proof}

Denote
$$
E_n(F)_{p,W}^*:=\inf_{T\in \mathcal{T}_n^*}\|F-T\|_{p,W}^*,
$$
where $\mathcal{T}_n^*$ is the set of all even trigonometric polynomials of degree at most $n$.

We need the following standard approximation inequalities.

\begin{lemma}\label{lemm2} (See~\cite[Theorems 4.1 and~4.3]{VIN2023}.)
Let $W\in A_p(\T)$, $1<p<\infty$,
$r\in\mathbb N$, and $\gamma>0$.
Then, for every even $F\in L_{p,W}(\T)$
and $n\in\mathbb N$, we have
\begin{equation}\label{eq:direct-symmetric}
E_{n-1}(F)_{p,W}^*
\le C\|(\I-A_{\gamma/n})^rF\|_{p,W}^*
\end{equation}
and
\begin{equation}\label{eq:inverse-symmetric}
\|(\I-A_{\gamma/n})^rF\|_{p,W}^*
\le
\frac{C}{n^{2r}}
\sum_{k=0}^{n-1}
(k+1)^{2r-1}E_k(F)_{p,W}^*.
\end{equation}
Here $C$ is independent of $F$ and $n$.
\end{lemma}

\begin{lemma}\label{lemKT}
Let $1<p<\infty$, $W\in A_p(\T)$,
$r\in\N$, and let $F\in L_{p,W}(\T)$
be even. Let $h>0$ and $n\in\N$
satisfy $n\asymp 1/h$, and suppose that
$T_n\in\mathcal T_{n-1}^*$ satisfies
\[
\|F-T_n\|_{p,W}^*
\le C E_{n-1}(F)_{p,W}^*.
\]
Then
\[
\|(\I-A_h)^rF\|_{p,W}^*
\asymp
\|F-T_n\|_{p,W}^*
+h^{2r}\|T_n^{(2r)}\|_{p,W}^*,
\]
where the constants are independent of
$F$, $n$, and $h$.
\end{lemma}

\begin{proof}
The assertion follows from
\cite[Theorem~4.6]{VIN2023}
combined with the weighted Bernstein
inequality for polynomials of near-best
approximation.
\end{proof}

\subsection{Relations for various $K$-functionals and moduli of smoothness}\label{secMK}

Let $P_{n,f}=P_n(f)_{p,(a,b)}\in\mathcal P_n$ be a polynomial of
near-best approximation to a function $f\in L_{p,(a,b)}$, that is,
$$
\|f-P_{n,f}\|_{p,(a,b)}
\leq
C E_n(f)_{p,(a,b)},
$$
where
$$
E_n(f)_{p,(a,b)}
:=
\inf_{P\in\mathcal P_n}
\|f-P\|_{p,(a,b)},
$$
and the constant $C$ is independent of $f$ and $n$.

\begin{lemma}\label{lemDTKM}
Let $1<p<\infty$, $a,b\ge0$, and
$r\in\mathbb N$.
Then, for every $f\in L_{p,(a,b)}$
and $t\in(0,t_0)$, we have
\begin{equation}\label{KMP}
K_{\varphi}^r(f,t)_{p,(a,b)}
\asymp \omega_{\varphi}^r(f,t)_{p,(a,b)}
\end{equation}
and
\begin{equation}\label{KRP}
\omega_{\varphi}^r(f,t)_{p,(a,b)}
\asymp
\|f-P_{n,f}\|_{p,(a,b)}
+t^r\|\varphi^rP_{n,f}^{(r)}\|_{p,(a,b)},
\qquad [1/t]=n.
\end{equation}
\end{lemma}

\begin{proof}
The first relation follows from~\cite[Theorem~6.1.1]{DITO1987}; see also~\cite[Theorem 2.5.1]{MAMI2008}. The second one is proved, e.g., in~\cite{GLTI2015}.
\end{proof}

\begin{lemma}\label{lemGTK}
(See~\cite[Theorem~6.2]{DI1998}.)
Let $1<p<\infty$, $r\in\N$, and suppose
that $a$ and $b$ satisfy~\eqref{condGT}.
Then, for every $f\in L_{p,(a,b)}$, we have
\[
{\rm K}_r(f,t)_{p,(a,b)}
\asymp
\|f-P_{n,f}\|_{p,(a,b)}
+t^r\|\mathcal D_rP_{n,f}\|_{p,(a,b)},
\qquad [1/t]=n.
\]
\end{lemma}

\begin{lemma}\label{lemGT}
(See~\cite[Corollary~2]{GLTI2015}.)
Let $1<p<\infty$, $r\in\N$, and suppose
that $a$ and $b$ satisfy~\eqref{condGT}.
Then, for every $f\in L_{p,(a,b)}$
and $t\in(0,t_0)$, we have
\[
C^{-1}K_\varphi^r(f,t)_{p,(a,b)}
\le {\rm K}_r(f,t)_{p,(a,b)}
\le C\left(
K_\varphi^r(f,t)_{p,(a,b)}
+t^r\|f\|_{p,(a,b)}
\right).
\]
\end{lemma}

\section{Proofs of the main results}

\subsection{Proof of Theorem~\ref{thsp}}

Put $h=\gamma/n$. By~\eqref{eqFf}--\eqref{eqFf+++}, it is sufficient to prove
\begin{equation}\label{en1}
\begin{split}
\|F-G_nF\|_{p,W}^*
\le C\bigl(
&\|(\I-A_h)^rF\|_{p,W}^*+\|(\I-A_h)^rF\|_{\ell_{p,W}(\Theta_n)}
\bigr).
\end{split}
\end{equation}
For simplicity, we write
\begin{equation}\label{not1}
\|\cdot\|_{p,W}:=\|\cdot\|_{p,W}^*,
\qquad
\|\cdot\|_{\ell_{p,W}}
:=\|\cdot\|_{\ell_{p,W}(\Theta_n)}.
\end{equation}
Define
\begin{equation}\label{not2}
F_h:=\sum_{k=1}^r
(-1)^{k+1}\binom rk A_h^kF
=F-(\I-A_h)^rF.
\end{equation}

By Lemma~\ref{lemmzF},
\begin{equation}\label{mo1}
\begin{split}
\|F-G_nF\|_{p,W}
&\le \|F-F_h\|_{p,W}
+\|F_h-G_nF_h\|_{p,W}+\|G_n(F-F_h)\|_{p,W}\\
&\le \|F-F_h\|_{p,W}
+\|F_h-G_nF_h\|_{p,W}+C\|F-F_h\|_{\ell_{p,W}}.
\end{split}
\end{equation}

Let $T\in\mathcal T_{n-1}^*$ be an even
trigonometric polynomial satisfying
\[
\|F-T\|_{p,W}
\le 2E_{n-1}(F)_{p,W}^*,
\]
and define
\[
T_h:=T-(\I-A_h)^rT.
\]
Since $A_h$ preserves even trigonometric
polynomials and their degrees, we have
$T_h\in\mathcal T_{n-1}^*$.
Consequently, $G_nT_h=T_h$.
Therefore, using Lemma~\ref{lemmzF},
\begin{equation}\label{zv1x}
\begin{split}
\|F_h-G_nF_h\|_{p,W}
&\le \|F_h-T_h\|_{p,W}
+C\|F_h-T_h\|_{\ell_{p,W}}.
\end{split}
\end{equation}

By the boundedness of the
averaging operators in Lemmas~\ref{lemsF} and~\ref{lemsFd},
\begin{equation}\label{zv2}
\|F_h-T_h\|_{p,W}
+\|F_h-T_h\|_{\ell_{p,W}}
\le C\|F-T\|_{p,W}.
\end{equation}
Indeed,
\[
F_h-T_h
=\sum_{k=1}^r
(-1)^{k+1}\binom rk A_h^k(F-T),
\]
and each term is bounded by successive
applications of the two lemmas.

Inequalities~\eqref{zv1x}
and~\eqref{zv2}, together with~\eqref{eq:direct-symmetric}, yield
\begin{equation}\label{mo2}
\|F_h-G_nF_h\|_{p,W}
\le C\|(\I-A_h)^rF\|_{p,W}.
\end{equation}
By~\eqref{not2},
\begin{equation}\label{mo2+}
\|F-F_h\|_{p,W}
=\|(\I-A_h)^rF\|_{p,W}.
\end{equation}
Combining~\eqref{mo1}, \eqref{mo2},
and~\eqref{mo2+}, we obtain~\eqref{en1},
which completes the proof.

\subsection{Proof of Theorem~\ref{prinvp}.}

As in the proof of Theorem~\ref{thsp},
we put $h=\gamma/n$ and use the
notation~\eqref{not1} and~\eqref{not2}.

By Lemma~\ref{lemmzF} and~\eqref{mo2},
\[
\begin{aligned}
\|F-F_h\|_{\ell_{p,W}}
&\le C\|G_n(F-F_h)\|_{p,W}\\
&\le C\bigl(
\|F-G_nF\|_{p,W}
+\|F-F_h\|_{p,W}
+\|F_h-G_nF_h\|_{p,W}
\bigr)\\
&\le C\bigl(
\|F-G_nF\|_{p,W}
+\|F-F_h\|_{p,W}
\bigr).
\end{aligned}
\]
Since $F-F_h=(\I-A_h)^rF$,
this yields~\eqref{st1invp}.

Furthermore,
\[
\begin{aligned}
\|(\I-A_h)^rF\|_{p,W}
&+\|(\I-A_h)^rF\|_{\ell_{p,W}}\\
&\quad\le C\bigl(
\|F-G_nF\|_{p,W}
+\|(\I-A_h)^rF\|_{p,W}
\bigr).
\end{aligned}
\]
By~\eqref{eq:inverse-symmetric},
\[
\|(\I-A_h)^rF\|_{p,W}
\le
\frac{C}{n^{2r}}
\sum_{k=0}^{n-1}
(k+1)^{2r-1}E_k(F)_{p,W}^*.
\]
Since $G_{k+1}F$ is an even
trigonometric polynomial of degree
at most $k$, we have
\[
E_k(F)_{p,W}^*
\le \|F-G_{k+1}F\|_{p,W}.
\]
Consequently,
\[
\begin{aligned}
&\|(\I-A_h)^rF\|_{p,W}
+\|(\I-A_h)^rF\|_{\ell_{p,W}}\\
&\quad\le C\left(
\|F-G_nF\|_{p,W}
+\frac1{n^{2r}}
\sum_{k=1}^n k^{2r-1}
\|F-G_kF\|_{p,W}
\right).
\end{aligned}
\]
By~\eqref{eqFf}--\eqref{eqFf+++},
this proves~\eqref{inv0p}.

\subsection{Proof of Theorem~\ref{thK1}}
In view of equalities~\eqref{eqFf}--\eqref{eqFf+++}, it is enough to prove that
\begin{equation}\label{eqKT1}
  \|F-K_n^r F\|_{p,W}^*\asymp\|(\I-A_{\g/n})^r F\|_{p,W}^*,
\end{equation}
where
$$
K_n^r F(t)=G_n F_h(t),\quad h=\g/n,
$$
and $F_h$ is defined by~\eqref{not2}.

In what follows, we use the notation~\eqref{not1}, as in the proof of Theorem~\ref{thsp}.
Applying inequality~\eqref{mo2}, we have
\begin{equation*}
\begin{split}
    \|F-K_n^rF\|_{p,W}&\le \|F-F_h\|_{p,W}+\|F_h-G_n F_h\|_{p,W}\\
    &\le C\|F-F_h\|_{p,W},
\end{split}
\end{equation*}
that proves the upper estimate in~\eqref{eqKT1}.

Further, let $T\in\mathcal T_{n-1}^*$ satisfy
\[
\|F-T\|_{p,W}\le 2E_{n-1}(F)_{p,W}^*.
\]
By~\eqref{zv1x} and~\eqref{zv2}, and since
$K_n^rF\in\mathcal T_{n-1}^*$, we obtain
\[
\begin{aligned}
\|F-F_h\|_{p,W}
&\le \|F-K_n^rF\|_{p,W}
+\|G_nF_h-F_h\|_{p,W}\\
&\le \|F-K_n^rF\|_{p,W}
+C\|F-T\|_{p,W}\\
&\le C\|F-K_n^rF\|_{p,W}.
\end{aligned}
\]
This proves the lower estimate in~\eqref{eqKT1}.

\subsection{Proof of Theorem~\ref{thbestapp}}

Denote $T(t)=P_n(\cos t)$ and $h=\gamma/n$.
Let $T_h$ be defined as in~\eqref{not2},
with $F$ replaced by $T$.
Using~\eqref{eqFf+}, \eqref{eqdFf},
and the notation~\eqref{not1}, we obtain
\begin{equation}\label{eqomega1}
\begin{split}
\|(\I-\A_h)^r f\|_{\ell_{p,w}(X_n)}
&=\|F-F_h\|_{\ell_{p,W}}\\
&\le \|F-T\|_{\ell_{p,W}}
+\|T-T_h\|_{\ell_{p,W}}+\|T_h-F_h\|_{\ell_{p,W}}\\
&=:I_1+I_2+I_3.
\end{split}
\end{equation}
By~\eqref{eqdFf}, we have
\begin{equation}\label{eqomega1+}
I_1=\|f-P_n\|_{\ell_{p,w}(X_n)}.
\end{equation}
By Lemmas~\ref{lemmzF} and~\ref{lemKT},
we obtain
\begin{equation}\label{eqomega2}
\begin{split}
I_2
&\le C\|T-T_h\|_{p,W}=C\|(\I-A_h)^rT\|_{p,W}\\
&\le Ch^{2r}\|T^{(2r)}\|_{p,W}\le C\|(\I-A_h)^rF\|_{p,W}.
\end{split}
\end{equation}
On the other hand, by Lemmas~\ref{lemsF}
and~\ref{lemsFd}, and the direct estimate
\eqref{eq:direct-symmetric}, we obtain
\begin{equation}\label{eqomega3}
\begin{split}
I_3
&\le C\|T-F\|_{p,W}\le C\|(\I-A_h)^rF\|_{p,W}.
\end{split}
\end{equation}

Combining~\eqref{eqomega1}--\eqref{eqomega3}
and using~\eqref{eqFf}, we obtain
\eqref{eqomega0}.

\subsection{Proof of Theorem~\ref{thOmega}}

Put
\[
J_\nu(f,n):=
\int_0^{1/n}
\frac{\Omega_\varphi^\nu(f,t)_{p,w}}
{t^{1+1/p}}\,dt.
\]
Let $M=\max\{n,\nu\}$ and
$M_j=2^jM$, $j\ge0$.
For every $j\ge0$, choose
$P_{M_j}\in\mathcal P_{M_j}$ satisfying
\[
\|f-P_{M_j}\|_{p,w}
\le 2E_{M_j}(f)_{p,w}.
\]

Since $w^{1/p}\in J_p^*$,
the integral Jackson estimate
\cite[Theorem~8.2.1]{DITO1987} gives
\begin{equation}\label{eq:omega-jackson}
\|f-P_{M_j}\|_{p,w}
\le C\int_0^{1/M_j}
\frac{\Omega_\varphi^\nu(f,t)_{p,w}}{t}\,dt.
\end{equation}

Set
\[
Q_j=P_{M_{j+1}}-P_{M_j}.
\]
Since $M_{j+1}\ge n$,
Lemma~\ref{lem:high-degree-sampling}
and~\eqref{eq:omega-jackson} yield
\[
\begin{split}
\|Q_j\|_{\ell_{p,w}(X_n)}
&\le C\left(\frac{M_{j+1}}n\right)^{1/p}
\|Q_j\|_{p,w}\\
&\le C\left(\frac{M_j}n\right)^{1/p}
\int_0^{1/M_j}
\frac{\Omega_\varphi^\nu(f,t)_{p,w}}{t}\,dt.
\end{split}
\]
Summing over $j$ and changing the
order of summation and integration,
we obtain
\[
\begin{split}
\sum_{j=0}^\infty
\|Q_j\|_{\ell_{p,w}(X_n)}
&\le Cn^{-1/p}
\int_0^{1/M}
\frac{\Omega_\varphi^\nu(f,t)_{p,w}}{t}
\sum_{\substack{j\ge0\\M_j\le1/t}}
M_j^{1/p}\,dt\\
&\le Cn^{-1/p}J_\nu(f,n).
\end{split}
\]

We next show that $P_{M_j}\to f$ locally uniformly
on $I$. By the weighted large-sieve inequality used
in the proof of Lemma~\ref{lem:high-degree-sampling},
we have
\[
|P(\cos\theta)|^p W_{2m}(\theta)
\le Cm\|P\|_{p,w}^p,
\qquad P\in\mathcal P_m.
\]
Since $W_{2m}$ is uniformly bounded below
on compact subintervals of $(0,\pi)$,
for every compact interval $K\subset(-1,1)$
we obtain
\[
\|P\|_{C(K)}
\le C_K m^{1/p}\|P\|_{p,w}.
\]
Applying this estimate to $Q_j$ and using
the preceding bounds, we obtain
\[
\sum_{j=0}^\infty\|Q_j\|_{C(K)}
\le C_K J_\nu(f,n)<\infty.
\]
Thus $P_{M_j}$ converges locally uniformly on $(-1,1)$.
Since $P_{M_j}\to f$ in $L_{p,w}(I)$ and
$f\in C_{\rm loc}(I)$, the limit must be $f$.
Consequently,
\[
f=P_M+\sum_{j=0}^\infty Q_j
\quad\text{on }(-1,1).
\]
In particular, the identity holds at
every interpolation node. Hence
\begin{equation}\label{eqomega8-}
\|f-P_M\|_{\ell_{p,w}(X_n)}
\le Cn^{-1/p}J_\nu(f,n).
\end{equation}
Also, by~\eqref{eq:omega-jackson},
\begin{equation}\label{eq:omega-remainder}
\|f-P_M\|_{p,w}
\le Cn^{-1/p}J_\nu(f,n).
\end{equation}

It remains to estimate
$\|(\I-\A_h)^r f\|_{\ell_{p,w}(X_n)}$.
We have
\begin{equation}\label{eq:omega-split}
\begin{aligned}
\|(\I-\A_h)^r f\|_{\ell_{p,w}(X_n)}
\le {}&
\|(\I-\A_h)^r P_M\|_{\ell_{p,w}(X_n)}\\
&+\|(\I-\A_h)^r(f-P_M)\|_{\ell_{p,w}(X_n)},
\end{aligned}
\end{equation}
where $h=\gamma/n$.
Since $(\I-\A_h)^r$
preserves $\mathcal P_M$ and $M/n\le\nu$,
Lemmas~\ref{lem:high-degree-sampling}
and~\ref{lemsF} yield
\begin{equation}\label{eq:omega-poly}
\begin{split}
\|(\I-\A_h)^rP_M\|_{\ell_{p,w}(X_n)}
&\le C\|(\I-\A_h)^rP_M\|_{p,w}\\
&\le C\bigl(
\|(\I-\A_h)^rf\|_{p,w}
+\|f-P_M\|_{p,w}
\bigr).
\end{split}
\end{equation}
On the other hand, applying Lemma~\ref{lemsFd}
to $\A_h^{j-1}(f-P_M)$ and then using
Lemma~\ref{lemsF}, we obtain
\[
\|\A_h^j(f-P_M)\|_{\ell_{p,w}(X_n)}
\le C\|f-P_M\|_{p,w},
\qquad j=1,\ldots,r.
\]
Hence, by the binomial expansion,
\begin{equation}\label{eq:omega-residual}
\begin{split}
&\|(\I-\A_h)^r(f-P_M)\|_{\ell_{p,w}(X_n)}
\le C\bigl(
\|f-P_M\|_{\ell_{p,w}(X_n)}
+\|f-P_M\|_{p,w}
\bigr).
\end{split}
\end{equation}
Combining~\eqref{eqomega8-}--\eqref{eq:omega-residual},
we obtain~\eqref{intOmega}.

\subsection{Proof of Corollary~\ref{corthsp}}

Let $P_{n-1}\in\mathcal P_{n-1}$ be a polynomial
of best approximation to $f$ in $L_{p,w}(I)$.
Since $\mathcal L_n(P_{n-1},u)=P_{n-1}$,
Lemma~\ref{lemmzf} yields
\[
\begin{aligned}
\|f-\mathcal L_n(f,u)\|_{p,w}
&\le \|f-P_{n-1}\|_{p,w}
+\|\mathcal L_n(f-P_{n-1},u)\|_{p,w}\\
&\le C\left(
\|f-P_{n-1}\|_{p,w}
+\|f-P_{n-1}\|_{\ell_{p,w}(X_n)}
\right).
\end{aligned}
\]
By the integral Jackson inequality
\cite[Theorem~8.2.1]{DITO1987}
and the telescoping argument used in the proof
of Theorem~\ref{thOmega}, now starting with
$M_j=2^j(n-1)$, we obtain
\[
\begin{aligned}
&\|f-P_{n-1}\|_{p,w}
+\|f-P_{n-1}\|_{\ell_{p,w}(X_n)}\le Cn^{-1/p}
\int_0^{1/(n-1)}
\frac{\Omega_\varphi^\nu(f,t)_{p,w}}
{t^{1+1/p}}\,dt.
\end{aligned}
\]
Finally, the fixed-factor scaling property
of the weighted main-part modulus,
which follows from the equivalence in
\cite[Lemma~8.1.2(E)]{DITO1987}
and the scaling argument in
\cite[Theorem~6.2.1]{DITO1987},
yields~\eqref{eq:lagrange-jackson}.

\subsection{{Proof of Theorem~\ref{th2-open-problem2}}}\label{secth2-open-problem2}

We first consider $n\ge2$.
Since $w$ is a Jacobi weight,
estimates~\eqref{eqmato1}
and~\eqref{eq:jacobi-local-weight},
together with the endpoint separation
in~\eqref{eqshi2}, yield
\[
\lambda_n(w,x_k)\asymp
\Delta_n(x_k)w(x_k),
\qquad k=1,\ldots,n.
\]
Moreover, by condition~\eqref{eqshi2}
and~\eqref{eqshi3},
for $k=2,\ldots,n-1$ and
$x\in[x_{k-1},x_{k+1}]$, we have
\[
\varphi(x)\asymp\varphi(x_k),
\qquad w(x)\asymp w(x_k),
\qquad
\Delta_n(x_k)\asymp\frac{\varphi(x_k)}{n}.
\]
Consequently,
\begin{equation}\label{newpr2}
\lambda_n(w,x_k)
\le \frac{C}{n}\varphi(x)w(x),
\quad x\in[x_{k-1},x_{k+1}],
\quad k=2,\ldots,n-1.
\end{equation}

Similarly, the endpoint spacing estimates
in~\eqref{eqshi2} and the Christoffel
estimates for Jacobi weights imply that,
for some sufficiently small $c>0$,
independent of $n$,
\begin{equation}\label{newpr3}
\lambda_n(w,x_1)
\le \frac{C}{n}\varphi(x)w(x),
\quad x\in[-1+cn^{-2},x_2],
\end{equation}
and
\begin{equation}\label{newpr4}
\lambda_n(w,x_n)
\le \frac{C}{n}\varphi(x)w(x),
\quad x\in[x_{n-1},1-cn^{-2}].
\end{equation}

Furthermore, by the mean value theorem, there exist points $\chi_{kn} \in [\theta_{kn}-\gamma/n, \theta_{kn}+\gamma/n] \subset [\theta_{k+1,n}, \theta_{k-1,n}]$ such that
$$
A_{\gamma/n}F(\theta_{kn}) = F(\chi_{kn}),
$$
which immediately implies
\begin{equation}\label{eqmeanv}
    \A_{\gamma/n} f(x_{kn}) = A_{\gamma/n}F(\theta_{kn}) = F(\chi_{kn}) = f(\xi_{kn}),
\end{equation}
where $\xi_{kn} = \cos(\chi_{kn}) \in  [\cos(\theta_{kn}+\gamma/n), \cos(\theta_{kn}-\gamma/n)]\subset [x_{k-1,n}, x_{k+1,n}]$ for $k=1,\ldots,n$.

By the superadditivity of the $p$-variation
and the continuity of $f$ in the interior of $I$,
for $-1<a<b<1$ we have
\begin{equation}\label{eq:variation-increment}
|f(b)-f(a)|^p
\le V_p(f,[a,b])
\le v_p(f,b)-v_p(f,a)
=\int_{(a,b]}dv_p(f,x).
\end{equation}
Applying Lemma~\ref{lemour} and~\eqref{eqmeanv}, we get
\begin{equation}\label{fffff}
    \begin{split}
        &\|f-\A_{\gamma/n} f\|_{\ell_{p,w}(X_n)}^p\\
        &\le C\sum_{k=1}^n \lambda_n(w,x_{k})|f(x_{k})- f(\xi_{k})|^p\\
        &\le C\bigg(\lambda_n(w,x_{1})\int_{\cos(\theta_{1n}+\gamma/n)}^{x_2}dv_p(f,x)+\sum_{k=2}^{n-1}\lambda_n(w,x_{k})\int_{x_{k-1}}^{x_{k+1}}dv_p(f,x)\\     &\qquad\qquad\qquad\qquad\qquad\qquad\qquad\qquad\,\,\,+\lambda_n(w,x_{n})\int^{\cos(\theta_{nn}-\gamma/n)}_{x_{n-1}}dv_p(f,x)\bigg)\\
        &=C(S_1+S_2+S_3).
    \end{split}
\end{equation}

By~\eqref{newpr2} and the bounded overlap
of the intervals $[x_{k-1},x_{k+1}]$, we obtain
\begin{equation}\label{newpr5}
S_2\le \frac{C}{n}
\sum_{k=2}^{n-1}
\int_{x_{k-1}}^{x_{k+1}}
\varphi(x)w(x)\,dv_p(f,x)
\le \frac{C}{n}
\int_{x_1}^{x_n}
\varphi(x)w(x)\,dv_p(f,x).
\end{equation}

To estimate $S_1$ and $S_3$, we use the well-known inequality
$$
\cos z\le 1-\frac{2}{\pi^2}z^2,\quad z\in [0,\pi],
$$
and Lemma~\ref{lemshi1},  which yield
\begin{equation*}
    \cos\({\theta_{1n}}+\frac\gamma n\)\ge \cos\(\pi-\frac{c_1-\g}{n}\){=}-\cos\(\frac{c_1-\g}{n}\)\ge -1+\frac{2}{\pi^2}\(\frac{c_1-\g}{n}\)^2
\end{equation*}
and
\begin{equation*}
    \cos\({\theta_{nn}}-\frac\gamma n\)\le \cos\(\frac{c_1-\g}{n}\)\le 1-\frac{2}{\pi^2}\(\frac{c_1-\g}{n}\)^2,
\end{equation*}
where, in view of~\eqref{gprob2}, the lower
constant $c_1$ in~\eqref{eqshi2}
is chosen so that $c_1>\gamma$.
 By~\eqref{gprob2}, we may choose the
constant $c>0$ in
\eqref{newpr3}--\eqref{newpr4}
sufficiently small so that
\[
0<c\le \frac{2}{\pi^2}(c_1-\gamma)^2.
\]
Consequently, the preceding estimates show
that the intervals of averaging at the
extreme nodes are contained in
$[-1+cn^{-2},1-cn^{-2}]$.
Applying~\eqref{newpr3} and~\eqref{newpr4},
we obtain
\begin{equation}\label{newpr6}
    S_1\le \frac{C}{n}\int_{-1+c/n^2}^{x_2}\vp(x)w(x)dv_{p}(f,x)
\end{equation}
and
\begin{equation}\label{newpr7}
    S_3\le \frac{C}{n}\int^{1-c/n^2}_{x_{n-1}}\vp(x)w(x)dv_{p}(f,x).
\end{equation}

Finally, combining~\eqref{fffff}, \eqref{newpr5}, \eqref{newpr6}, and \eqref{newpr7}, we obtain
\begin{equation*}
       \|f-\A_{\gamma/n} f\|_{\ell_{p,w}(X_n)} \le \frac{C}{n^{1/p}} \left( \int_{-1+c/n^2}^{1-c/n^2} w(x)\varphi(x) \, dv_p(f,x) \right)^{1/p},
\end{equation*}
which completes the proof for $n\ge 2$.
For $n=1$, the assertion follows directly
from the same variation estimate applied
to the single node.

\subsection{{Proof of Theorem~\ref{th2-open-problem2+}}}\label{secth2-open-problem2+}

Put $h=\gamma/n$ and $F(t)=f(\cos t)$.
By Jensen's inequality,
\[
|F(\theta)-A_hF(\theta)|^p
\le \frac1{2h}
\int_{\theta-h}^{\theta+h}
|f(\cos\theta)-f(\cos v)|^p\,dv.
\]

By~\eqref{eq:variation-increment}, for almost
every $\theta$ and $v$, we have
\[
|f(\cos\theta)-f(\cos v)|^p
\le \int_{J(\theta,v)}dv_p(f,x),
\]
where $J(\theta,v)$ is the interval with
endpoints $\cos\theta$ and $\cos v$.
Writing
\[
\int_{J(\theta,v)}dv_p(f,x)
=
\int_{-1}^1
\mathbf 1_{\{x\in J(\theta,v)\}}\,dv_p(f,x)
\]
and applying Fubini's theorem, we obtain

\begin{equation}\label{eeqq1}
\begin{aligned}
\|f-\mathcal A_hf\|_{p,w}^p
&\le \frac1{2h}
\int_{-1}^1\int_0^\pi
w(\cos\theta)\sin\theta\\
&\quad\times
\int_{\theta-h}^{\theta+h}
\mathbf 1_{\{x\in J(\theta,v)\}}
\,dv\,d\theta\,dv_p(f,x).
\end{aligned}
\end{equation}

For $v\in\mathbb R$, let $\rho(v)\in[0,\pi]$
be determined by $\cos\rho(v)=\cos v$.
For every $\theta\in[0,\pi]$,
\[
|\rho(v)-\theta|\le |v-\theta|.
\]
Fix $x=\cos t$, $t\in(0,\pi)$.
If $x\in J(\theta,v)$, then $t$ lies between
$\theta$ and $\rho(v)$, and hence
\[
|\theta-t|\le |v-\theta|\le h.
\]
Since the inner integral with respect to $v$
is at most $2h$, it follows from~\eqref{eeqq1}
that
\begin{equation}\label{eeqq2}
\|f-\mathcal A_hf\|_{p,w}^p
\le \int_{-1}^1
\left(
\int_{\max(0,t-h)}^{\min(\pi,t+h)}
w(\cos\theta)\sin\theta\,d\theta
\right)dv_p(f,x),
\end{equation}
where $x=\cos t$.

By Taylor's formula, for $x=\cos t$ and
$|\theta-t|\le h=\gamma/n$, we have
\[
|\cos\theta-x|
\le \sin t\,|\theta-t|
+\frac{|\theta-t|^2}{2}
\le \frac{\gamma}{n}\varphi(x)
+\frac{\gamma^2}{2n^2}
\le C_\gamma\Delta_n(x).
\]
Setting
\[
a=\max(0,t-h),\qquad
b=\min(\pi,t+h),
\]
and making the change of variables
$y=\cos\theta$, we obtain
\[
\begin{aligned}
\int_a^b w(\cos\theta)\sin\theta\,d\theta
&=\int_{\cos b}^{\cos a}w(y)\,dy\\
&\le
\int_{I\cap[x-C_\gamma\Delta_n(x),
             x+C_\gamma\Delta_n(x)]}
w(y)\,dy.
\end{aligned}
\]
Since $w$ is a doubling weight, the last
integral is bounded by a constant multiple
of the integral over the interval with
radius $\Delta_n(x)$. Consequently,
\begin{equation}\label{eeqq3}
\begin{aligned}
&\int_{\max(0,t-h)}^{\min(\pi,t+h)}
w(\cos\theta)\sin\theta\,d\theta\le C
\int_{I\cap[x-\Delta_n(x),x+\Delta_n(x)]}
w(y)\,dy.
\end{aligned}
\end{equation}

Combining~\eqref{eeqq2} and~\eqref{eeqq3}
gives~\eqref{eqeq1} and completes the proof.

\subsection{Proof of Theorem~\ref{thArK}}

Let $F(t)=f(\cos t)$ and
$W^{(a,b)}(t)=w^{(a,b)}(\cos t)|\sin t|.$
The assumptions on $a$ and $b$ imply that
$W^{(a,b)}\in A_p(\T)$.

Let $m=[1/h]$ and
$T(t)=P_{m,f}(\cos t)$, where $P_{m,f}$
is a polynomial of near-best approximation.
Applying Lemma~\ref{lemKT} with $n=m+1$,
we obtain
\[
\begin{aligned}
\|(\I-\A_h)^r f\|_{p,(a,b)}
&=\|(\I-A_h)^rF\|_{p,W^{(a,b)}}^*\\
&\asymp
\|F-T\|_{p,W^{(a,b)}}^*
+h^{2r}\|T^{(2r)}\|_{p,W^{(a,b)}}^*.
\end{aligned}
\]
Since
\[
T^{(2r)}(t)
=(-1)^r(\mathcal D_{2r}P_{m,f})(\cos t),
\]
it follows that
\[
\begin{aligned}
\|(\I-\A_h)^r f\|_{p,(a,b)}
\asymp\,
&\|f-P_{m,f}\|_{p,(a,b)}+h^{2r}\|\mathcal D_{2r}P_{m,f}\|_{p,(a,b)}.
\end{aligned}
\]
The assertion now follows from
Lemma~\ref{lemGTK} applied with order $2r$.

%%%%%%%

\section*{Acknowledgments}

{\small The authors (L. Angeloni and D. Costarelli) are members of the Gruppo Nazionale per l'Analisi Matematica, la Probabilit\`a e le loro Applicazioni (GNAMPA) of the Istituto Nazionale di Alta Matematica (INdAM), of the network RITA (Research ITalian network on Approximation), and of the UMI (Unione Matematica Italiana) group T.A.A. (Teoria dell'Approssimazione e Applicazioni).
}

\section*{Funding}

{\small The authors (L. Angeloni and D. Costarelli) have been partially supported within the project "National Innovation Ecosystem grant ECS00000041 - VITALITY", funded by the European Union - NextGenerationEU under the Italian Ministry of University and Research (MUR).
}

{\small Yu. Kolomoitsev is supported by the German Research Foundation, project KO 5804/3-1.}

\section*{Conflict of interest/Competing interests}

{\small The authors declare that they have no competing interests.}

\section*{Availability of data and material and Code availability}

{ \small Not applicable.}

%

%%%%%%

\begin{thebibliography}{99}

\bibitem{ABM} J. Appell, J. Bana\'s and N. Merentes, Bounded variation and around, De Gruyter Series in Nonlinear Analysis and Applications, De Gruyter, Berlin,  17, 2014.

\bibitem{AS1973} R. Askey, Summability of Jacobi series, Trans. Amer. Math. Soc. \textbf{179} (1973), 71--84.

\bibitem{BA1974} V. M. Badkov, Convergence in the mean and almost everywhere of Fourier series in polynomials orthogonal on an interval, Math. USSR-Sb. \textbf{24} (1974), No.~2, 223--256.

\bibitem{BU1983} P. L. Butzer, Legendre transform methods in the solution of basic problems in algebraic approximation, Functions, series, operators I, II (Budapest, 1980) 277--301; Colloq. Math. Soc. Janos Bolyai, vol. 35, North-Holland, Amsterdam, 1983.

\bibitem{CHDI1994} W. Chen, Z. Ditzian, Strong converse inequality for Kantorovich polynomials,
Constr. Approx. \textbf{10} (1994), No.~1, 95--106.

\bibitem{CHDIIV1993} W. Chen, Z. Ditzian, K. Ivanov, Strong converse inequality for the Bernstein-Durrmeyer operator, J. Approx. Theory \textbf{75} (1993), No.~1, 25--43.


\bibitem{DI1998}
Z. Ditzian, Fractional derivatives and best approximation, Acta Math. Hungar. \textbf{81} (1998), No.~4, 323--348.

\bibitem{DI2007} Z. Ditzian,
Polynomial approximation and $\omega_\vp^r(f,t)$ twenty years later,
Surv. Approx. Theory \textbf{3} (2007), 106--151.

\bibitem{DITO1987} Z. Ditzian, V. Totik, Moduli of smoothness, Springer-Verlag, New York (1987).

\bibitem{GLTI2015} P. Glazyrina, S. Tikhonov,
Jacobi weights, fractional integration, and sharp Ulyanov inequalities, J. Approx. Theory 195 (2015), 122--140.

\bibitem{GOZH1995} H. Gonska, X.-l. Zhou, The strong converse inequality for Bernstein-Kantorovich operators, Comput. Math. Appl. \textbf{30} (1995), 103--128.


\bibitem{IV1982} K. G. Ivanov,  On a new characterization of functions I, Serdica, \textbf{8} (1982), 262--279.

\bibitem{IV1986} K. G. Ivanov,  Converse theorems for approximation by Bernstein polynomials in $L_p[0,1]$ $(1<p<\infty)$,
Constr. Approx. \textbf{2} (1986), 377--392.


\bibitem{JOU2017}
H.~Joung,
Weighted inequalities for generalized
polynomials with doubling weights,
J. Inequal. Appl. (2017), Article 91.


\bibitem{K25}  Yu. Kolomoitsev,  Special measures of smoothness for approximation by sampling operators in $L_p(\R^d)$, J. Math. Anal. Appl. \textbf{559} (2026), No.~1, 130439.


\bibitem{KO2026} Y.~Kolomoitsev,
Approximation by linear sampling operators in Banach spaces,
J. Fourier Anal. Appl.
\textbf{32} (2026), Article 93.

\bibitem{KL23}	Yu. Kolomoitsev, T. Lomako,  Sharp  $L_p$-error estimates for sampling operators, J. Approx. Theory \textbf{294} (2023), 105941.


\bibitem{KLP2018}
Y. Kolomoitsev, T. Lomako, J. Prestin,
On $L_p$-error of bivariate polynomial
interpolation on the square,
J. Approx. Theory \textbf{229} (2018), 13--35.


\bibitem{KOTI2026} Yu. Kolomoitsev, S. Tikhonov,
Marcinkiewicz-Zygmund inequalities in quasi-Banach function spaces,
J. Complexity \textbf{93} (2026),  101999.

\bibitem{MAMI2008} G. Mastroianni, G. Milovanovi\'c, Interpolation Processes. Basic Theory and Applications, Springer, Berlin,
2008.


\bibitem{MARU1999}
G. Mastroianni, M. G. Russo,
Lagrange interpolation in weighted Besov spaces,
Constr. Approx. \textbf{15} (1999), No.~2, 257--289.

    
\bibitem{MARU2010}
G. Mastroianni, M. G. Russo,
Some new results on Lagrange interpolation for bounded variation functions,
J. Approx. Theory \textbf{162} (2010), 1417--1428.

\bibitem{MATH2008} G. Mastroianni, W. Themistoclakis,  De la Vall\'ee Poussin means and Jackson's theorem, Acta Sci. Math. \textbf{74} (2008), No.~1-2, 147--170.

 \bibitem{MATO2000}
G. Mastroianni, V. Totik,
Weighted polynomial inequalities with doubling and $A_\infty$ weights,
Constr. Approx. \textbf{16} (2000), No.~1, 37--71.

\bibitem{MATO2010}
G.~Mastroianni and V.~Totik,
Uniform spacing of zeros of orthogonal polynomials,
Constr. Approx. \textbf{32} (2010), 181--192.

\bibitem{NEVAI1979}
P. Nevai, Orthogonal polynomials,
Mem. Amer. Math. Soc. \textbf{18} (1979), No.~213.


\bibitem{NE1984} P. Nevai, Mean convergence of Lagrange interpolation. III, Trans. Amer. Math. Soc. \textbf{282} (1984), No.~2, 669--698.

\bibitem{POKA1998} M. K. Potapov, G. N. Kazimirov, On the approximation of functions that have a given order of $k$th generalized modulus of smoothness by algebraic polynomials, Mat. Zametki \textbf{63} (1998), no. 3, 425--436; translation in Math. Notes \textbf{63} (1998), no. 3--4, 374--383.

\bibitem{PR93} J. Prestin,  Lagrange interpolation for functions of bounded variation,
Acta Math. Hung. \textbf{62} (1993), No.~1-2, 1--13.


\bibitem{SHI2005}
Y. G. Shi, Christoffel type functions for $m$-orthogonal polynomials,
J. Approx. Theory \textbf{137} (2005), No.~1, 57--88.

\bibitem{ST1993} E. M. Stein. Harmonic analysis: real-variable methods, orthogonality and oscillatory integrals.
Princeton: Princeton Univ. Press. 1993.

\bibitem{TH2018} W. Themistoclakis, Weighted  $L_1$ approximation on  $[-1,1]$ via discrete de la Vall\'ee Poussin means, Math. Comput. Simul. \textbf{147} (2018), 279--292.

   
\bibitem{VIN2023}
O. L. Vinogradov,
Direct and inverse theorems of approximation theory
in Lebesgue spaces with Muckenhoupt weights,
Ufim. Mat. Zh. \textbf{15} (2023), No.~4, 42--60;
translation in Ufa Math. J. \textbf{15} (2023), No.~4, 42--61.

\end{thebibliography}
\end{document}